\documentclass[11pt]{article}

\usepackage[latin1]{inputenc}
\usepackage[dvips]{graphicx}
\usepackage{latexsym}
\usepackage{amsmath,amssymb,amsfonts,amsthm}
\usepackage{xcolor}

\usepackage{latexsym}
\usepackage{amssymb}
\usepackage[dvips]{graphicx}

\newtheorem{Theorem}{Theorem}[part]
\newtheorem{Definition}{Definition}[part]
\newtheorem{Proposition}{Proposition}[part]

\newtheorem{Lemma}{Lemma}[part]

\def \cqfd{\quad_\diamondsuit}

\def \R{\mathbb{R}}

\def \E{\mathbb{E}}

\def \P{\mathbb{P}}
\def \Q{\mathbb{Q}}

\def \ni{\noindent}

\title{Convergence of multi-dimensional quantized $SDE$'s}

\date{}

\begin{document}

\author{Gilles PAG{\`E}S~\thanks{LPMA, Universit\'e Paris 6, case 188,  4 pl. Jussieu, 75252 Paris Cedex 5, France
\texttt{gilles.pages@upmc.fr}} \and
Afef SELLAMI~\thanks{JP Morgan, London \& LPMA, Universit\'e Paris 6, case 188,  4 pl. Jussieu, 75252 Paris Cedex 5, France \texttt{afef.x.sellami@jpmorgan.com}}}

%
%
\maketitle

\begin{abstract}We quantize a multidimensional  
$SDE$ (in the Stratonovich sense) by solving the related system of $ODE$'s in which the $d$-dimensional Brownian motion has been replaced by the components of functional stationary quantizers. We
make a connection  with rough path theory to show that the solutions of the quantized solutions of the $ODE$  converge toward the solution of the $SDE$. On our way to this result we provide  convergence rates of optimal quantizations toward the Brownian motion for $\frac 1q$-H\" older distance, $q>2$, in $L^p(\P)$. 
\end{abstract}



\vspace{7mm}


\vspace{5mm}


\section{Introduction }

Quantization  is a way to discretize the path space of a random phenomenon: a random vector in finite dimension, a stochastic
process in infinite dimension.  Optimal Vector Quantization theory  
(finite-dimensional)  random vectors  finds its origin in the
early 1950's in order to discretize some emitted signal (see~\cite{GRLU}). It was further developed by specialists in Signal Processing and later in
Information Theory.  The infinite dimensional case started to be  extensively investigated in the
early 2000's by several authors (see $e.g.$  \cite{LUPA1}, \cite{DEFEMASC}, \cite{LUPA2}, \cite{LUPA3}, \cite{LUPA2}, \cite{DER},  \cite{GRLUPA2}, etc).

In  \cite{LUPA3}, the functional quantization of a class
of Brownian diffusions has been investigated from a constructive point of view.
The main feature of this class of diffusions was that the diffusion coefficient was the
inverse of the gradient of a diffeomorphism (both coefficients being
smooth). This class contains most (non degenerate) scalar diffusions.
Starting from a sequence of rate optimal quantizers, some
sequences of quantizers of the Brownian diffusion are produced as solutions of (non coupled)
$ODE$'s. This approach relied on the Lamperti transform and was closely
related to the Doss-Sussman representation formula of the flow of a
diffusion as a functional of the Brownian motion. In many situations
these quantizers are rate optimal (or almost rate optimal) $i.e.$ that
they quantize the diffusion at the same rate $O((\log N)^{-\frac 12})$ as
the Brownian motion itself  where $N$ denotes the generic size of the
quantizer. In a companion paper (see~\cite{PAPR2}),   some
cubature formulas based on some of  these quantizers were implemented, namely those
obtained from some optimal product quantizers based on the
Karhunen-Lo\`eve expansion of the Brownian motion, to price some Asian
options in a Heston stochastic volatility model. Rather unexpectedly in
view of the theoretical rate of convergence, the numerical experiments
provided quite good numerical results for some ``small" sizes of quantizers.
Note however that these numerical implementations  included some
further speeding up procedures combining the stationarity of the quantizers and the Romberg
extrapolation leading to a $O((\log N)^{-\frac 32})$ rate.  Although this result relies on some still pending conjectures about the asymptotics of bilinear functionals of the quantizers, it strongly pleads in favour of the construction of such stationary (rate optimal) quantizers, at least  when one has in mind possible numerical applications. 

Recently  a sharp quantization rate  ($i.e.$ including an explicit constant) has been established for a class of  not too degenerate $1$-dimensio\-nal Brownian diffusions. However the approach is not   constructive   (see~\cite{DER}). On the other hand, the
standard rate $O((\log N)^{-\frac 12})$ has been extended in~\cite{LUPA5}
to general
$d$-dimensional  It\^o processes, so including $d$-dimensional Brownian diffusions regardless of their ellipticity properties. This latter approach, based   an expansion in the Haar basis,  is
constructive, but the resulting quantizers are no longer stationary.

Our aim in this paper is to extend the constructive natural approach initiated
in~\cite{LUPA3} to general $d$-dimensional diffusions in order to produce some rate optimal stationary quantizers of these processes. To this end, we will
 call upon some seminal results from  rough path theory, namely the continuity of the It\^o map,  to replace the ``Doss-Sussman setting".  In fact we will show that if one replaces in an $SDE$ (written in the Stratonovich sense) the Brownian motion by some elementary quantizers, the solutions of the resulting $ODE$'s make up  some rough paths  which converge (in $p$-variation and in the H\" older metric) to the solution of the $SDE$. We use her the rough path theory as a tool and we do not aim at providing new insights on this theory. We can only mention that these rate optimal stationary quantizers can be seen as  a new example of rough paths, somewhat less ``stuck" to a true path of the underlying process. 
 
 This work is devoted to Brownian diffusions which is naturally the prominent example in view of applications, but it seems clear that this could be extended to $SDE$ driven $e.g.$ by fractional Brownian motions (however our approach requires to have an explicit form for the Karhunen-Lo\`eve basis as far as numerical implementation is concerned).

 \medskip Now let us be more precise.  We consider a diffusion process
\[
dX_t = b(t,X_t)\,dt+\sigma(t,X_t)\circ dW_t, \; X_0= x\!\in\R^d,\;
t\in [0,T],
\]
in the Stratonovich sense where $b: [0,T]\times \R^d\to \R^d$ and  $\sigma:[0,T]\times
\R^d\to {\cal M}(d\times d)$ are continuously differentiable with linear growth (uniformly with respect to
$t$) and $W=(W_t)_{t\in [0,T]}$ is a $d$-dimensional Brownian motion defined on a filtered probability space $(\Omega,{\cal A}, \P)$. (The fact that the state space and $W$ have the same dimension is in no case a restriction since our result has nothing to do with ellipticity).

Such an $SDE$ admits a unique strong solution denoted
$X^x =(X^x_t)_{t\in [0,T]}$ (the dependency in $^x$ will be dropped from now to alleviate notations). The $\R^d$-valued  process $X$ is pathwise continuous  and
 $\sup_{t\in [0,T]}|X_t|\!\in L^r(\P)$, $r>0$ (where $|\,.\,|$
denotes the canonical Euclidean norm on $\R^d$). In particular $X$ is  
bi-measurable and can be seen as an $L^r(\P)$-Radon random variable taking
values in the Banach spaces $(L^p_{T,\,\R^d},|\,.\,|_{L^p_{_T}})$ where
$L^p_{T,\,\R^d}=L^p_{\R^d}([0,T], dt)$ and $|g|_{L^p_{_T}}=
\left(\int_0^T |g(t)|^pdt\right)^{\frac 1p}$ denotes the usual
$L^p$-norm when $p\!\in[1,\infty)$. 

\smallskip For every integer $N\ge 1$, we can investigate for $X$ the level $N$
$(L^r(\P),L^p_{_T})$-quantization problem for this process $X$, namely solving the
minimization of the $L^r(\P)$-mean $L^p_{T,\R^d}$-quantization error
\begin{equation}\label{OptiQuant}
e_{N,r}(X,L^p):=\min\left\{e_{N,r}(\alpha,X,L^p),\; \alpha\subset L^p_{T,
\R^d},\,\mbox{card}\,\alpha\le N\right\}
\end{equation}
where $e_{N,r}(\alpha,X,L^p)$ denotes the $L^r$-mean quantization error induced by $\alpha$, namely
\[
e_{N,r}(\alpha,X,L^p):= \left(\E
\min_{a\in \alpha}|X-a|^r_{p}\right)^{\frac 1r} =
\left\|\min_{a\in \alpha}|X-a|_{L^p_{T,\R^d}}\right\|_{L^r(\P)}.
\]

The use of ``$\min$" in~(\ref{OptiQuant}) is justified by the existence of an optimal quantizer solution to that problem as shown in~\cite{CUMA, GRLUPA3} in this infinite dimensional setting.
The Voronoi diagram associated to a quantizer $\alpha$ is a Borel partition $(C_a(\alpha))_{a\in \alpha}$ such that
\[
C_a(\alpha)\subset \Big\{x\!\in L^p_{T,\R^d}\,|\, |x-a|_{L^p_{T,\R^d}}\le \min_{b\in \alpha}|x-b|_{L^p_{T,\R^d}}\Big\}
\]
and a functional quantization of $X$ by $\alpha$ is defined by the nearest neighbour projection of $X$ onto $\alpha$ related to the Voronoi diagram
\[
\widehat X^{\alpha}:= \sum_{a\in \alpha} a \mbox{\bf 1}_{\{X\in C_a(\alpha)\}}.
\]

In finite dimension (when considering $\R^d$-valued random vectors   instead of $L^p_{T,\R^d}$-valued 
processes) the answer is provided   by the so-called
Zador Theorem which says (see~\cite{GRLU}) that if
$\E|X|^{r+\delta}<+\infty$ for some $\delta>0$ and if $g$ denotes the
absolutely continuous part of its distribution  then
\begin{equation}\label{Zador}
N^{\frac 1d} e_{N,r}(X,\R^d) \to \widetilde J_{r,d}
\|g\|^{\frac 1r}_{\frac{d}{d+r}}\qquad\mbox{as} \qquad N\to \infty
\end{equation}
where $ \widetilde J_{r,d}$  is finite positive real constant obtained as the limit of the normalized
quantization error when $X\stackrel{d}{=}U([0,1])$. This constant is
unknown except when $d=1$ or $d=2$.

A  non-asymptotic version of Zador's Theorem can be found $e.g.$  in~\cite{LUPA5}: for every $r,\, \delta>0$
there exists a universal constant $C_{r,\delta}>0$ and an integer $N_{r,\delta}\ge $ such that, for every random vector $\Omega,{\cal A}, \P)\to \R^d$,  
\[
\forall\, N\ge N_{r,\delta},\qquad  e_{N,r}(X,\R^d) \le 
C_{r,\delta} \|X\|_{r+\delta} N^{-\frac 1d}.
\]
The asymptotic behaviour of the $L^s(P)$-quantization error of sequences of $L^r$-optimal quantizers of a random vector $X$ when $s>r$  has been extensively investigated in~\cite{GRLUPA3} and will be one crucial tool to establish our mains results.

In infinite dimension, the case of Gaussian processes was the first to
have been extensively investigated, first in the  purely quadratic case
($r=p=2$):   sharp rates have been
established for a wide family of Gaussian processes including the
Brownian motion,   the fractional Brownian motions 
(see~\cite{LUPA1,LUPA2}). For these two processes   sharp rates are
also known for  $p\!\in[1,\infty]$ and
$r\!\in (0,\infty)$ (see~\cite{DER}).  More recently, a connection
between mean   regularity of $t\mapsto X_t$ (from $[0,T]$ into
$L^r(\P)$) and the  quantization rate has been established (see~\cite{LUPA5}): if the above mapping is
$\mu$-H\"older for an index $\mu\!\in (0,1]$, then 
\[
e_{N,r}(X,L^p) = O((\log N)^{-\mu}),\qquad p\!\in (0,r).
\]
Based on this result, some universal quantization rates have been obtained for general L\'evy
processes with or without Brownian component some of them turning out to be optimal, once compared with the
lower bound estimates derived from small deviation theory (see~$e.g.$~\cite{GRLUPA1} or~\cite{DEFEMASC}).  One
important feature of interest of the purely quadratic case  is that it is possible to construct from the
Karhunen-Lo\`eve expansion of the process  two families of rate optimal (stationary) quantizers, relying on 

\smallskip
--  sequences $(\alpha^{(N, prod)})_{N\ge 1}$ of optimal {\em product
quantizers} which are rate optimal $i.e.$ such that 
$e_{N,r}(\alpha^{(N)},X,L^2)= O(e_{N,2}(X,L^2))$ (although not with a sharp optimal rate).

\smallskip
--  sequences of {\em true } optimal quantizers  (or at least some good
numerical approximations) $(\alpha^{(N,*)})_{N\ge 1}$ $i.e.$ such that
$e_{N,r}(\alpha^{(N,*)},X,L^2)= e_{N,2}(X,L^2)$.

\bigskip
We refer to Section~\ref{Bckgrnd} below for further insight on these objects
(both being available on the website {\tt www.quantize.math-fi.com}).

\medskip
The main objective of this paper is the following: let
$(\alpha^{N})_{N\ge 1}$ denote a sequence of rate optimal stationary (see~(\ref{stationnaire}) further on) quadratic quantizers of
a $d'$-dimensional standard Brownian motion $W=(W^1,\ldots,W^{d})$.
Define the sequence  $ x^{N}=(x^{N}_n)_{ n=1,\ldots,N}$, $N\ge 1$, of solutions of the ODE's
\[
x^{N}_n(t)= x+\int_0^t b(x^{N}_n(s))ds +\int_0^t
\sigma(x^{N}_n(s))d\alpha_n^{N}(s),\quad n=1,\ldots,N.
\]

Then, the finitely valued-process defined by 
\[
\widetilde X^N =\sum_{n=1}^N x^{N}_n \mbox{\bf 1}_{\{W\!\in
C_n(\alpha^{(N)})\}}
\]
converges toward the diffusion $X$ on $[0,T]$ (at least  in probability) as $N\to
\infty$.  This convergence will hold with respect to distance introduced  in the rough path theory  (see~\cite{LYO, LEJ, FRI, FRIVIC, LYO3})  which always implies  convergence with respect to the $\sup$ norm. The reason is that our result will  appear as an application of  (variants of the) the celebrated {\em Universal Limit Theorem} originally established by T. Lyons in~\cite{LYO}.  The distances of interest in rough path theory
are related to  the  $\frac 1q$-H\"older semi-norm   or the $q$-variation semi-norm both when $q>2$ defined for every $x\!\in {\cal C}([0,T],\R^d)$ by
\[
\|x\|_{q,Hol}= T^{\frac 1p}\sup_{0\le s< t\le T}\frac{|x(t)-x(s)|}{|t-s|^{\frac 1q}}\le +\infty,
\]
and 
\begin{eqnarray*}
 {\rm Var}_{q, [0,T]}(x)&\!:=\! &
\sup\Big\{\Big(\sum_{0\le \ell\le k-1}|x(t_{\ell+1})\!-\!   x(t_\ell)|^q\Big)^{\!\frac
1q}\\
&&\hskip 2 cm  ,0\le t_0\le t_1\le \dots \le t_k \le T,k\ge 1\!\Big\} \le +\infty
\end{eqnarray*} 
respectively. Note that 
\[
\|x-x(0)\|_{\rm \sup} \le  {\rm Var}_{p, [0,T]}(x)\le \|x\|_{p,Hol}.
\]
 
 From a technical viewpoint we aim at applying some continuity results  established on the It\^o map by several authors (see $e.g.$~\cite{LEJ, LYO,FRI, LEJ2}) that is the continuity of a solution $x$ of the $ODE$ (in a rough path sense)
 \[
 dx_t = f(x_t)dy_t, \quad x_0=x(0),
 \]
 as a functional of $y$. However, the above (semi-)norms associated to a function $x$ are not sufficient  and the natural space to define such  rough ODE is not the  ``naive" space of  paths but a space of enhanced paths, which involves in the case of a multi-dimensional  Brownian motion  the mutual L\'evy areas of its components.  Convergence in this space is defined by considering appropriate $\frac 1q$-H\" older and $p$-variation semi-norms to  both the $d$-dimensional path and the related (pseudo-)L\'evy areas (wit different values of $q$ and $p$, see Section~\ref{CvgQSDE}).   Our application to  quantized $SDE$'s   will make extensively use the fact that our functional quantizations of the Brownian motion $W$ will all satisfy a stationary assumption $i.e.$
 \[
 \widehat W = \E (W\,|\,\sigma (\widehat W))
 \]
so that we will extend the Kolmogorov criterion satisfied by $W$ to its functional quantizers $\widehat W$ for free. This approach is rather straightforward and its field of application seems more general than our functional quantization purpose: thus the piecewise affine interpolations of the Brownian motion obviously satisfy such a property (see Appendix).

\medskip The paper is organized as follows. In Section~2 we provide some short background on functional quantization as well as  preliminary elementary results on stochastic integration with respect to a stationary functional quantizer of a $d$-dimensional standard Brownian motion. In Section~3, we define a quantized approximation scheme of an $SDE$ (in the Stratonovich sense) driven by a standard Brownian motion by its functionally quantized counterpart which turns out to be a system of (non-coupled) $ODE$'s. To this end we recall some basic facts on rough path theory, in particular the notion of convergence we need to define on the so-called {\em multiplicative functionals} involved in the continuity of the It\^o map which, when dealing with Brownian motion amounts, to some convergence in H\" older semi-norm of the naive path as well as, roughly speaking,  the running (pseudo-)L\'evy areas of its components.  In Section~4 and~5, we establish successively the convergence in the H\" older distance of sequences of optimal stationary  quantizations $\widehat W$ of the Brownian motion toward $W$: Section~4 is devoted to the  convergence of the ``regular"  paths whereas  Section~5  deals with the convergence  of the running (pseudo-)L\'evy areas (and to the global convergence of the couple). In both cases we provide some convergence rate in the $(\log N)^{-a}$, $a\!\in (0, \frac 12)$ scale which is the natural scale for such convergences since optimal functional quantizations  of the Brownian motion are known to converge at a $(\log N)^{-\frac 12}$-rate  for most usual norms (like quadratic pathwise norm on  $L^2([0,T],dt)$).

\bigskip
\noindent{\sc Notations:} $\bullet$ For every $d\ge 1$, one denotes  $\xi = (\xi^1,\ldots,\xi^d)$ a row vector of $\R^d$. ${\cal M}(d\times d)$ will denote the set of square matrices with $d$ lines.

\noindent $\bullet$ $|\,.\,|$ denotes the canonical
Euclidean norm on
$\R^d$.

\noindent $\bullet$ We denote $({\cal F}^X_t)_{t\ge 0}$ the augmented natural filtration of a process $X=(X_t)_{t\ge 0}$  (so that it satisfies 
the usual conditions).  

\noindent $\bullet$ For a bounded function $f:[0,T]\to \R^d$,  $\|\,f\,\|_{\rm sup}:= \sup_{t\in [0,T]}|f(t)|$. If $f$ is a Borel function and $p\!\in[1,+\infty)$,  $\|f\|_{L^p_{T,\R^d}}:=\left(\int_0^T|f(t)|^pdt\right)^{\frac 1p}$.

\noindent $\bullet$ For an $\R^d$-valued bi-measurable process $X$ and   $p\!\in [1,+\infty)$, we  denote  $\|X\|_{p}:=\|\,|X|_{L^p_{T,\R^d}}\,\|_{p}=\left(\E\int_0^T |X_t|^pdt\right)^{1/p}$.

\noindent $\bullet$ We  denote  $t^n_k =\frac{kT}{2^n}$, $k=0,\ldots,2^n$, the uniform mesh of the interval $[0,T]$, $T>0$ and  $I^n_k =[t^n_k,t^n_{k+1}]$, $k=0,\ldots,2^n-1$.  

\noindent $\bullet$ $\lfloor x\rfloor$ denotes the lower integral part of $x\!\in \R$.

\noindent $\bullet$ Let $(a_n)_{n\ge 0}$ and $(b_n)_{n\ge 0}$ be two sequences of real numbers: $a_n  \sim b_n $ if  $a_n = b_n+o(b_n)$ and $a_n \asymp b_n $  if  $a_n =O(b_n)$ and $b_n= O(a_n)$.

\section{Background and preliminary results on functional quantization}
\subsection{Some background on functional quantization}\label{Bckgrnd}
Functional quantization of stochastic processes can be seen as a discretization of the path-space of a process and the
 approximation (or coding) of a process by finitely many deterministic functions from its path-space. In a Hilbert
space setting this reads as follows.

Let $(H, \langle \cdot , \cdot \rangle) $ be a separable Hilbert space with norm
$|  \cdot | $
and let $X : (\Omega, {\cal A}, \P) \rightarrow H$ be a random vector taking its values in $H$ with distribution $\P_X$.
Assume the integrability condition
\begin{equation}\label{(1.2)}
\E \,|  X |^2 < +\infty.
\end{equation} 
For $N\! \ge 1$, the $L^2$-optimal $N$-quantization problem for $X$  
consists in minimizing
\[
\Big \|\min_{a \in \alpha} | X-a |  \Big\|_{L^2(\P)} = \left(\E \min_{a \in \alpha} | X - a |^2 \right)^{1/2} 
\]
over all subsets $\alpha \subset H$ with $\mbox{card} (\alpha) \leq N$. Such a set $\alpha$ is called
$N$-codebook or $N$-quantizer.  The minimal  quantization
error of $X$ at level $N$ is then defined by
\begin{equation}\label{infimum}
e_{_N}(X,H) := \inf \left\{ ( \E \min_{a \in \alpha} | X-a |^2)^{1/2} : 
\alpha \subset H, \; \mbox{card} (\alpha)
\leq N \right\} .
\end{equation}
For a given $N$-quantizer $\alpha$ one defines an associated nearest neighbour projection
\[
\pi_\alpha := \sum\limits_{a \in \alpha} a \mbox{\bf 1}_{C_{a} (\alpha)}
\]
and the induced $\alpha$-(Voronoi)quantization   of $X$ by setting 
\begin{equation}
\hat{X}^\alpha := \pi_\alpha (X) ,
\end{equation}
where $\{ C_a (\alpha) : a \in \alpha\}$ is a Voronoi partition induced by $\alpha$, that is a Borel partition of $H$ satisfying
\begin{equation}
C_a (\alpha) \subset 
\{ x \in H : | x-a | = \min_{b \in \alpha} | x-b |  \}
\end{equation}
for every $a \in \alpha$. Then one easily checks that, for any random vector $X^{'} : \Omega \rightarrow \alpha \subset H$,
\[
\E\,| X - X^{'} |^2 \geq \,\E\,| X - \hat {X}^\alpha |^2 = \E\,\min_{a \in \alpha} | X-a |^2
\]
so that finally
\begin{eqnarray}
\hskip -0.35 cm e_n (X,H) \!&=&\!\inf \left\{\!  \Big\| |X - q(X)| \Big\|_{L^2(\P)}\! ,\, q : H \stackrel{Borel}{\rightarrow} H, \mbox{card} (q(H)) \leq N \right\} \nonumber \\
        \!&=&\!\inf \left\{ \! \Big\|   |X- Y | \Big\|_{L^2(\P)}\!,\,Y\!:\! (\Omega,{\cal A}) \stackrel{r.v.}{\rightarrow} H,
             \mbox{card} (Y(\Omega)) \leq N \right\}\!.  \;\;
\end{eqnarray}

A typical setting for functional quantization is $H\!=\! L^2_{_T}\!:=\!L_{\R}^2([0,1],\!dt)$ (equipped with $\langle f,g\rangle_2:=\int_0^Tfg(t)dt$ and     $|f|_{L^2_T} :=
\sqrt{\langle f,f\rangle_2}$). Thus any (bi-measurable, real-valued) process 
$X=(X_t)_{t\in[0,T]}$    defined on a probability space  $(\Omega,{\cal A}, \P)$ such that 
\[
\int_0^T \E(X^2_t) dt <+\infty 
\]
is a random variable $X: (\Omega,{\cal A}, \P)\to L^2_{_T}$. But this Hilbert setting is not the only possible
one for functional quantization (see $e.g.$~\cite{LUPA4},~\cite{GRLUPA2},~\cite{DEFEMASC}, etc) since natural
Banach spaces like $L_{\R}^p([0,T],dt)$ or ${\cal C}([0,T], \R)$ are natural path-spaces. 

In  the purely Hilbert setting   the existence of (at least) one {\em optimal
$N$-quantizer} for every integer $N\ge 1$  is established so that  the infimum in~(\ref{infimum}) holds as a
minimum. A typical feature of this quadratic Hilbert framework is the so-called {\em stationarity} (or self-consistency) property  satisfied by such an optimal $N$-quantizer $\alpha^{(N,*)}$:
\begin{equation}\label{stationnaire}
  \hat X^{\alpha^{(N,*)}} = \E (X\,|\,   \hat X^{\alpha^{(N,*)}} ).
 \end{equation}
 This property, known as stationarity,  will be used extensively throughout the paper.

  This existence property holds true in any reflexive Banach space and $L^1$ path spaces (see~\cite{GRLUPA2} for details).
\subsection{Constructive aspects of functional quantization of the   Brownian motion}
  \subsubsection{Karhunen-Lo\`eve basis ($d= 1$)} 
First we consider a scalar Brownian motion $(W_t)_{t\in [0,T]}$ on a probability space $(\Omega,{\cal A}, \P)$. The two
main classes of rate optimal quantizers of the Brownian motion are  the product optimal quantizers and the true
optimal quantizers. Both are based on the Karhunen-Lo\`eve expansion of the Brownian motion given by
\begin{equation}\label{KL}
W_t = \sum_{k\ge 1}\sqrt{\lambda_k}\, \xi_k\, e^W_{k}(t)
\end{equation}
where, for every $k\ge 1$,  
\begin{equation}\label{KL2}
\lambda_k =
\left(\frac{T}{\pi(k-1/2)}\right)^2  \quad
\mbox{ and } \quad  e^W_{k}(t)= \sqrt{\frac{2}{T}}\sin\left(\frac{t}{\sqrt{\lambda_k}}\right)
\end{equation}
and  
$$
\xi_k= \frac{(W\,|\, e_k^W)_{_2}}{\sqrt{\lambda_k}} = \sqrt{\frac{2}{T}}\int_0^T W_t \,\sin(t/\sqrt{\lambda_k})\frac{dt}{\sqrt{\lambda_k}}.
$$
The sequence $(e^W_{k})_{k\ge 1}$ is an orthonormal basis of $L^2_{_T}$. The system $(\lambda_k,e^W_k)_{k\ge 1}$ can be
characterized as the eigensystem of the symmetric positive trace class covariance operator of $f\mapsto(t\mapsto  \int_0^T
(s\wedge t)\,f(s)ds)\equiv(t\mapsto   \E(<\!f\,|\,W\!>_2W_t)$. In particular this implies  that the
Gaussian sequence
$(\xi_k)_{k\ge 1}$ is  pairwise uncorrelated hence i.i.d.,  
${\cal N}(0;1)$-distributed.  The Karhunen-Lo\`eve expansion of $W$ plays the role of $PCA$ of the process: it is the fastest way to 
exhaust  the  variance of $W$ among all expansions on an orthonormal basis.

The convergence of the series in the right hand side of~(\ref{KL})  holds  in  $L^2_{_T}$ for every $\omega
\!\in\Omega$  and $\P(d\omega)$-$a.s.$ for every
$t\!\in [0,T]$. In fact this convergence also holds
in $L^2(\P)$ and $\P(d\omega)$-$a.s.$ for the sup norm over
$[0,T]$. The first convergence  follows from  Theorem~\ref{Thm1}$(a)$ further on applied with $X=W$ and ${\cal
G}_N = \sigma(\xi_1,\ldots,\xi_{_N})$ and the second one follows $e.g.$ from~\cite{LUPA4} $\P(d\omega)$-$a.s.$.   In
particular the convergence holds in $L^2(d\P\otimes dt)$ or equivalently in  $L^2_{L^2_{_T}}(\P)$.
Note that this basis has already been used in the framework of rough path theory for Gaussian processes, see $e.g.$~\cite{COVI, FRIVIC01,FRIVIC02}.

\subsubsection{Optimal product quantization ($d\ge 1$)} 

\noindent $\rhd$ {\em The one-dimensional case $d=1$.} The previous expansion of the Brownian motion suggests to define a product quantization of $W$ at level $N$ by
\begin{equation}\label{prodquantifB}
\widehat W^{(N_1,\ldots,N_{_L})}_t := \sqrt{\frac{2}{T}}\sum_{k= 1}^{L}\sqrt{\lambda_k}\, \widehat
\xi^{N_k}_k
\sin\left(\frac{t}{\sqrt{\lambda_k}}\right)
\end{equation} 
where $N_1,\, \ldots,N_{_L}$ are non zero integers satisfying $N_1\cdots N_{L}\le
N$ and $\widehat
\xi^{N_1}_1,\ldots,\widehat\xi_{_L}^{N_{_L}}$ are optimal quadratic quantizations of
$\xi_1,\ldots,\xi_{_L}$. The resulting (squared) quadratic quantization error reads
\begin{equation}\label{OptiProdQuant}
\|W-  \widehat W^{(N_1,\ldots,N_{_L})}\|^2_{_2}=\sum_{k\ge 1} \frac{\lambda_k}{N_k^2}.
\end{equation}
An {\em optimal product $N$-quantization} $\widehat W^{N,prod}$ is obtained as a
solution to the following {\em integral bit allocation optimization}  problem  for the sequence $(N_k)_{k\ge 1}$: 
\begin{equation}\label{BitAlloc}
 \min\left\{\!\|W\!-\!  \widehat W^{(N_1,\ldots,N_{_L})}\|_{_2},\,N_1, \ldots,N_{_L}\ge 1,\,N_1\!\cdots\!
N_{L}\le N,\,L\!\ge\!1 \!\right\}
\end{equation}
(see~\cite{LUPA1} for further details and~\cite{PAPR2} for the numerical aspects). It is established in~\cite{LUPA1} (as a special case  of a   more general
result on Gaussian processes) that 
\begin{equation}\label{prodquantrate}
\frac 1 T\|W- \widehat W^{N,prod}\|_{_2}  \asymp (\log N)^{-\frac 12}
\end{equation} 
Furthermore,   the critical dimension $L=L_W(N)$ satisfies $L_W(N) \sim
\log N$. Numerical experiments carried out in~\cite{PAPR2} show   that  
$$
\frac 1T \|W- \widehat
W^{N, prod}\|_{_2}\approx c_{_W} (\log N)^{-\frac 12}
$$ 
with $c_{_W}\approx 0.5$ (at least up to $N\le 10 000$). 

It is possible to get a closed form for  the underlying optimal product quantizers $\alpha^N$. First, note  that the normal distribution on the real line being
$\log$-concave, there is exactly  one stationary  quadratic quantizer of full size $M$ for every $M\ge 1$ (hence it is the optimal one). So, let $N\ge 1$ and
let
$(N_k)_{k\ge 1}$ denote its optimal integral bit  allocation for the Brownian motion $W$. For every $N_k\ge 1$, we denote by
$\beta^{(N_k)}:=\{\beta^{(N_k)}_{i_k},\;1\le i_k\le N_k\}$  the unique optimal quantizer of the normal distribution: thus $\alpha{(0)}=\{0\}$ by symmetry of
the normal distribution.   Then, the optimal quadratic product $N$-quantizer $\alpha^{N,prod}$ (of ``true size" $N_1\times\cdots \times N_{L_W(N)}\le N$) can
be described using a multi-indexation as follows:
\[
\alpha^{N,prod}_{(n_1,\ldots,n_k,\ldots)}(t)=\sum_{k\ge 1}\beta_{n_k}^{(N_k)}\sqrt{ \lambda_k} e^W_k(t),\qquad n_k\!\in\{1,\ldots,N_k\},\;\;k\ge 1.
\]
These sums are in fact all finite so that all the functions  $\alpha^{N,prod}_{(i_1,\ldots,i_n,\ldots)}$ are  ${\cal C}^\infty$ with finite variation on
every interval of $\R_+$.

Explicit optimal \textsl{integral  bit allocations}  as well as  optimal quadratic quantizations (quantizers and their weights) of the scalar normal distribution are
available on the website~\cite{website}. Note for practical applications that  this  optimal product quantization is based on $1$-dimensional quantizations of
small size of the scalar  normal distribution ${\cal N}(0;1)$. This kind of functional quantization  has been
 applied in~\cite{PAPR2} to   price Asian options in a Heston stochastic volatility model.

 \medskip 
 \noindent  $\rhd$ {\em The $d$-dimensional case.} Assume now $W= (W^1,\ldots,W^d)$ is a $d$-dimensional Brownian motion. Its optimal product quantization at level $N\ge 1$ will be defined as the optimal product quantization at level $\lfloor N^{\frac 1d}\rfloor$ of each of its $d$ components.
 
 \medskip
\noindent $\rhd$ {\em Additional results on optimal vector quantization of the normal distribution on $\R^d$.}   We will extensively make use of the  {\em distortion mismatch} result established in~\cite{GRLUPA3} that we recall here only  in the $d$-dimensional Gaussian case. Let $Z$ be an ${\cal N}(0;I_d)$ random vector and let $\alpha^N$ be an optimal quadratic quantizer at level $N$ of $Z$ (hence of size $N$). Then
 \begin{eqnarray}
 \label{match} (i)& \forall\, p\!\in (0,2+d),\,\forall\, N\ge 1,\,\|Z-\widehat Z^{\alpha^N}\|_p \le C_{Z,p}N^{-\frac 1d}, \hskip 2cm  \\ \nonumber \\
 \nonumber (ii) & \forall\, p\!\in [2+d,+\infty),\,\forall\, \eta\!\in (0,d+2),  \forall\, N\ge 1,\hskip 3cm\\
\label{mismatch} &\|Z-\widehat Z^{\alpha^N}\|_p \le C_{Z,p,\eta}N^{-\frac{2+d-\eta}{dp}} 
 \end{eqnarray}
where $C_{Z,p}$ and $C_{Z,p,\eta}$ are two positive real constants.
\subsubsection{Optimal  quantization ($d=1$)} It is established in~\cite{LUPA1} (Theorem~3.2)  that the quadratic
optimal quantization of the one-dimensional Brownian motion reads
\begin{equation}\label{quantifB}
\widehat W^{N,opt}_t = \sqrt{\frac{2}{T}}\sum_{k= 1}^{d_W(N)}\sqrt{\lambda_k}\, 
(\widehat \zeta_{d_W(N)}^N)^k \sin\left(\frac{t}{\sqrt{\lambda_k}}\right)
\end{equation}
where, for every integer $d\ge 1$,  $\zeta_{d} = {\rm Proj}_{E_{d}}^\perp\!(W)\!\sim \!{\cal
N}\!\left(0; {\rm Diag}(\lambda_1,\ldots,
\lambda_{d})\!\right)$ with $E_{d}:=\R$-${\rm span
}\left\{\sin\left(./\sqrt{\lambda_1}\right),\ldots,\sin\left(./\sqrt{\lambda_{d}}\right) \right\}$
and $\widehat \zeta_{d}^N$ is an optimal quadratic quantization of  $\zeta_{d}$ at level  (or of size) $N$. 

\medskip If one considers an optimal quadratic $N$-quantizer $\beta^N=\{\beta^{N}_n,\; n=1,\ldots,N\}\subset\R^{d_W(N)}$ of the distribution ${\cal
N}\left(0; {\rm Diag}(\lambda_1,\ldots,
\lambda_{d_W(N)})\right)$ ({\em a priori}  not unique) 
\[
\alpha^{N,opt}_{n}(t)=\sum_{k= 1}^{d_W(N)}(\beta^{(N)}_{n})^k\sqrt{ \lambda_k} \,e^W_k(t),\qquad n=1,\ldots,N.
\]
Once again this defines a ${\cal C}^\infty$ function with finite variation on
every interval of $\R_+$.

A sharp rate has been obtained in~\cite{LUPA2} for the  resulting optimal quantization error  
\begin{equation}\label{optiquantif}
  \|W-\widehat W^{N,opt}\|_{_2} \sim T c_{_W}^{opt}   (\log N)^{-\frac 12} \quad \mbox{ as }\quad N\to \infty
\end{equation}
where $c_{_W}^{opt}=\frac{\sqrt{2}}{\pi}  \approx 0.4502  $. 

The    true value of the critical dimension $d_W(N)$ is unknown. A 
conjecture supported by numerical evidences is that $d_W(N)\sim
\log N$. Recently a first step to this conjecture has been established in~\cite{LUPAWI} by showing that  
$$
\liminf_N\frac{d_W(N)}{\log(N)} \ge \frac 12. 
$$ 

\medskip
Large scale computations of   optimal quadratic 
quantizers of the Brownian motion have been carried out (up to $N= 10\,000$ and $d=10$). They  are 
available on the website~\cite{website}. 

In the $d$-dimensional setting, several definitions of an {\em optimal quantization} of the Brownian motion $W=(W^1,\ldots,W^d)$ can be given. For our purpose, it is  convenient to adopt  the following one:  
$$
\widehat W^{N,opt}:= \Big(  \widehat{W^i}^{\lfloor N^{\frac 1d}\rfloor, opt}\Big)_{1\le i\le d}.
$$
Its property of interest is that this definition preserves the component\-wise independence as well as a stationarity property (see below) since
\[
\E \Big(W^i \,|\, \widehat W^{N,opt}\Big) =\E \Big(W^i \,|\, \widehat{W^i}^{\lfloor N^{\frac 1d}\rfloor, opt}\Big) = \widehat{W^i}^{\lfloor N^{\frac 1d}\rfloor, opt}, \; i=1,\ldots,d.
\]

\subsubsection{Wiener like integral with respect to a stationary  functional  quantization ($d=1$)}\label{1.2.3}
Both types of quantizations  defined above share an important property of quantizers: stationarity.

\begin{Definition} Let $\alpha\subset L^2_{_T}$, $\alpha\neq \emptyset$, be a quantizer. The quantizer $\alpha$ is stationary for the  (one-dimensional) Brownian motion $W$ if there is a Voronoi quantization  $\widehat W:=\widehat W^\alpha$ induced by $\alpha$ such that
\begin{equation}\label{stationarity}
\widehat  W= \E(W\,|\, \sigma(\widehat W))\qquad a.s.
\end{equation}
 where $\E(\,.\,|{\cal G})$ denotes the functional conditional expectation given the
$\sigma$-field ${\cal G}$ on $L_{L^2_{_T}}^2(\P)$ (see Appendix) and $\sigma(\widehat W)$ is the $\sigma$-field spanned by $\widehat W$. 
\end{Definition}

Note that if $\alpha$ is stationary for one Brownian motion, so it is for any Brownian motion since this stationarity property only depends on the  Wiener distribution. 

\medskip
 In the case of product quantization $\widehat W^{N,prod} $, this follows from the stationarity
property of the optimal quadratic quantization of the marginals $\xi_n$ (see~\cite{LUPA1} or~\cite{PAPR2}). In the case of optimal quadratic quantization $\widehat W^{N,opt}$ this follows
from  the optimality of the quantization of    $\zeta_{d_W(N)}$ itself. 

\medskip
We will  now  define  a kind of {\em Wiener integral with respect} to such  {\em a stationary quantization} $\widehat W$ of a {\em one-dimensional} $W$. So we assume that  $d=1$ until the end of this Section.

\medskip
First, we must have in mind that if $W$ is an $({\cal F}_t)$-Brownian motion where the filtration $({\cal F}_t)_{t\ge 0}$ satisfies the usual conditions,  one can define the Wiener stochastic integral (on $[0,T]$) of any process $\varphi\!\in L^2([0,T]\times\Omega, {\cal B}([0,T])\otimes {\cal F}_0, dt \otimes d\P )$ with respect to $W$. The non-trivial case is when ${\cal F}^W_t \neq {\cal F}_t$, typically when ${\cal F}_t = {\cal F}^B_{T} \vee {\cal F}^W_t$, $t\!\in [0,T]$ where $B$ and $W$ are independent.  One can see it as a special case of It\^o stochastic integral or as an extended Wiener integral: if $(\varphi(t,\omega))_{(\omega,t)\in \Omega\times [0,T]}$ denotes an elementary process of  the form
\[
\varphi(t,\omega):= \sum_{k=1}^n \varphi_k(\omega)\mbox{\bf 1}_{s_k<t \le s_{k+1}},\; 0=s_0<s_1<\cdots<s_{n-1}<s_n=T
\]
where the random variables $\varphi_i$ are ${\cal F}_0$-measurable (hence independent of $W$). Set 
$$
I_{_T}(\varphi):=\displaystyle 
\sum_{k=1}^n \varphi_k(W_{s_{k+1}}-W_{s_k}).
$$ 
Then, $I_{_T}$ is an isometry  from $L^2_{L^2_{_T}}(\P)$ into $L^2({\cal F}_{_T},\P)$.
Furthermore, one easily checks that
\[
\E\left(\int_0^T \varphi(s,.)dW_s\,|\, {\cal F}^W_{_T}\right)= \int_0^T\E\left( \varphi(s,.)\,|\, {\cal F}^W_{_T}\right)dW_s
\] 
where ${\cal F}^W_{_T}$ denotes the augmented filtration of $W$ at time $T$. We follow the same lines to define the stochastic integral with respect to a stationary quantizer.  Set for the same  elementary process $\varphi$
\[
\widehat I_{_T}(\varphi) = \sum_{k=1}^n \xi_k (\widehat W_{s_{k+1}}-\widehat W_{s_k})
\]
so that
\begin{eqnarray*}
\widehat I_{_T}(\varphi) &=& \sum_{k=1}^n \xi_i\,\E (W_{s_{k+1}}-W_{s_k}\,|\,\widehat  W)\\
 &=& \sum_{k=1}^n\E\big( \xi_k (W_{s_{k+1}}-W_{s_k})\,|\,{\cal F}_0\vee\sigma(\widehat  W)\big)\\
&=&\E\left(\int_0^T\varphi(t,.) dW_t\,|\,{\cal F}_0\vee\sigma(\widehat  W) \right)
\end{eqnarray*}
where we used that  the $\sigma$-fields $\sigma(\widehat W)$ and ${\cal F}_0$ are independent since $\widehat W$ is a Borel function of $W$.  As a consequence,
\[
\|\widehat I_{_T}(\varphi)\|_{_2}^2\le \|I_{_T}(\varphi)\|_{_2}^2= \|\,|\varphi|_{L^2_{_T}}\|_{_2}^2.
\]
Hence, the linear transformation $\widehat I_{_T}$ extends into a linear continuous mapping on the whole set  $L^2_{L^2_{_T}}({\cal F}_0, \P)$. Furthermore, one checks, first on
 elementary processes, then on   $L^2_{L^2_{_T}}({\cal F}_0, \P)$ by continuity of the (functional)  conditional expectation, that
\[
\E\left( I_{_T}(\varphi)\,|\,{\cal F}_0\vee\sigma(\widehat  W) \right)
=\widehat I_{_T}(\varphi).
\] 
We will denote from now on $\widehat I_{_T}(\varphi)(\omega) $ as an integral, namely
\[
\widehat I_{_T}(\varphi)(\omega) := \int_0^T \varphi(t,\omega)d\widehat W_t(\omega).
\]

Now  set as usual, for every $t\!\in [0,T]$,  
$$
\displaystyle \int_0^t \varphi(s,\omega)d\widehat W_s(\omega):= \int_0^T\mbox{\bf 1}_{[0,t]}(s)\varphi(s,\omega)d\widehat W_s(\omega).
$$ 

One checks  using Jensen and Doob Inequality that,
\begin{eqnarray}
\nonumber \E\sup_{t\in [0,T]}\left|\int_0^t \varphi(s,.)d\widehat W_s\right|^2&\le &\E\sup_{t\in[0,T]}\left|\int_0^t\varphi(s,.)dW_s\right|^2\\
 &\le& 4 \,\E \int_0^T\varphi^2(s,.) \, ds.\label{JenDoob}
\end{eqnarray}
Furthermore, as soon as the underlying stationary quantizer $\alpha$ (such that $\widehat W= \widehat W^\alpha$) is made up with  pathwise continuous elements,  for every  elementary process  $\varphi$, its integral process 
\[
\int_0^t \varphi(s,.) \,dW_s= \sum_{k=1}^n \xi_k(\widehat W_{s_{k+1}\wedge t}-\widehat W_{s_k \wedge t})
\]
pathwise continuous as well since $\widehat W$ is $\alpha$-valued. One classically  derives, by combining this result with~(\ref{JenDoob}) and the everywhere density of elementary processes, that, for every $\varphi\!\in L^2_{L^2_{_T}}({\cal F}_0, \P)$, the process
\[
\left(\int_0^t \varphi(s,.)d\widehat W_s\right)_{t\in [0,T]}\qquad \mbox{ admits a continuous modification.}
\]
This is  always  this modification that will be considered from now on.  As a matter of fact, if $\varphi_n$ denotes a sequence of elementary processes in $L^2_{L^2_{_T}}({\cal F}_0, \P)$ converging to $\varphi$, $i.e.$ satisfying
\[
\E\,\int_0^T (\varphi-\varphi_n)^2(s,.) ds \longrightarrow 0\qquad \mbox{ as }\; n\to \infty.
\]
It follows from~(\ref{JenDoob}) that the convergence also holds in $L^2_{L^\infty_{_T}}({\cal F}_0, \P)$. In particular, there is a subsequence that converges $\P$-$a.s.$ for the $\|\,.\,\|_{\sup}$ which implies the existence of a continuous modification  for $\displaystyle \int_0^t \varphi(s,\omega)d\widehat W_s(\omega)$.

Finally, using the characterization of functional conditional expectation (see Appendix), it follows that
\begin{equation}\label{espcondfonc}
\E\left(  \int_0^. \varphi(s,.)d\widehat W_s  \ ,|\,{\cal F}_0\vee\sigma(\widehat  W) \right)
= \int_0^.\varphi(s,.)d\widehat W_s.
\end{equation}

\begin{Proposition}\label{Lem1.1} Let $W$ be a (real-valued) ${\cal F}_t$-standard Brownian motion. 

\medskip
\noindent $(a)$ For every $\varphi\!\in L^2_{L^2_{_T}}({\cal F}_0, \P)$
\begin{equation}\label{Expand1}
\int_0^t \varphi(s,.)dW_s= \sqrt{\frac 2T}\sum_{k\ge 1}\xi_k  \int_0^t \varphi(s,.) \cos(s/\sqrt{\lambda_k})ds
\end{equation}
where $\xi_k:=(W|e^W_k)_{_2}/\sqrt{\lambda_k}$ are independent, ${\cal N}(0;1)$-distributed (see~(\ref{KL}) and~(\ref{KL2})) and independent of $\varphi$.

\medskip
\noindent $(b)$ Let $\widehat W$ be a stationary quantization of $W$. For every $\varphi\!\in L^2_{L^2_{_T}}({\cal F}_0, \P)$
\begin{equation}\label{Expand2}
\int_0^t \varphi(s,.)d\widehat W_s= \sqrt{\frac 2T}\sum_{k\ge 1} \frac{(\widehat W|e^W_k)_{_2}}{\sqrt{\lambda_k}} \int_0^t \varphi(s,.)
\cos(s/\sqrt{\lambda_n})ds.
\end{equation}
\end{Proposition}
In particular if $\widehat W$ is a product quantization, then 
$$
\frac{(\widehat W|e^W_k)_{_2}}{\sqrt{\lambda_k}}=\widehat{\frac{( W|e^W_k)_{_2}}{\sqrt{\lambda_k}}}=\widehat \xi_k.
$$

\bigskip
\noindent {\bf Proof.}  $(a)$   Set for every $\varphi\!\in L^2_{L^2_{_T}}({\cal F}_0, \P)$, 
\begin{eqnarray}\label{Wienerintexp}
J_{_T}(\varphi)&:=& \sqrt{\frac 2T} \sum_{k\ge 1} \xi_k \sqrt{\lambda_k}\int_0^T \varphi(s,.)d\sin(s/\sqrt{\lambda_k})\\
\nonumber &=& \sqrt{\frac 2T} \sum_{k\ge 1} \xi_k \int_0^T \varphi(s,.)\cos (s/\sqrt{\lambda_k})ds.
\end{eqnarray}
 This defines clearly an isometry from $L^2_{L^2_{_T}}({\cal F}_0, \P)$ into the Gaussian space spanned by $(\xi_n)_{n\ge 1}$ since
\[
\E(J_{_T}(\varphi)^2) = \frac 2T \sum_{k\ge 1} \E(\xi_k^2) \,\E \left(\int_0^T\hskip -0.25cm g(s)\frac{1}{\sqrt\lambda_k}
\cos(s/\sqrt{\lambda_k})ds\right)^2 \hskip -0.15cm = \E \int_0^T\hskip -0.25cm g^2(t)dt.
\]
The last equality uses that the sequence  $\Big(\sqrt{\frac 2T}\cos(\pi(k-\frac 12)t/T)\Big)_{k\ge 1}$ is an orthonormal basis of $L^2_{_T}$. Finally, note that  for
every
$t\!\in [0,T]$, $J_{_T}(\mbox{\bf 1}_{[0,t]}) =\sqrt{\frac 2T}\sum_{k\ge 1}\sqrt{\lambda_k}\,\xi_k\sin(t/\sqrt{\lambda_k}) =
W_t$. This proves that   $J_{_T}=I_{_T}$ $i.e.$ is but the (extended) Wiener integral with respect to $W$.  

\medskip
\noindent $(b)$ This follows by taking the (functional) conditional expectation of~(\ref{Expand1}). $\cqfd$

 \subsubsection{Application  to multi-dimensional Brownian motions $(d\ge 2$)}

Now we apply the above result to a componentwise (stationary) functional quantization of a multi-dimensional standard Brownian motion.

\begin{Proposition}\label{Prop1.1} Let $W=: (W^1,\ldots, W^d)$ denote a $d$-dimensional standard Brownian motion and let $\widehat W:= (\widehat W^1,\ldots, \widehat W^d)$ be a pathwise continuous
stationary quantization of $W$ (no optimality is requested here).  Then, $\P$-$a.s.$, for every $i\neq j$, $i,\,j\!\in \{1,\ldots,d\}$, for every $s$,$ t\!\in [0,T]$, $0\le s\le t$, 
\[
\E\left(\int_s^t(W^i_u-W^i_s)dW^j_u\,|\, \sigma(\widehat W)\right) = \int_s^t(\widehat W^i_s-\widehat W^i_s) d\widehat W^j_u.
\]
\end{Proposition}

\noindent {\bf Proof.}  All the components of $\widehat W$ being independent, it is clear one can replace $\sigma(\widehat W)$ by $\sigma (\widehat W^i, \widehat W^j)$. Then, the stochastic integral $\int_0^. W^i_s dW^j_s$ coincides with the (extended) Wiener integral defined with respect 
to the filtration ${\cal G}^{j}_{i,t}:=\sigma({\cal F}^{W^i}_{_T},{\cal F}^{W^j}_t)$ (it is clear that $W^j$ is a ${\cal G}^{j}_{i,t}$-standard Brownian
motion still by independence). The result is then a straightforward consequence of~(\ref{espcondfonc}).$\cqfd$

\bigskip
\ni{\bf Remark.} The above result still holds if one considers an additional ``$0^{th}$"~component $W^0_t\!=\!t$  to the Brownian motion and to its functional quantization by setting $\widehat W^0_t\!=\!t$ as well.

\section{Convergence of quantized $SDE$'s: a rough path approach}\label{CvgQSDE}

\subsection{From It\^o to Stratonovich}
An $SDE$ 
\[
dX_t = b(t,X_t)dt +\sigma(t,X_t) dW_t, \quad X_0\!\in L^p_{\R^d}(\P)
\]
where $b:[0,T]\times \R^d\to \R^d$ and $\sigma:[0,T]\times \R^d \to {\cal M}(d\times q)$ are smooth
enough functions ($e.g.$ continuously differentiable with bounded differentials) and $W=(W_t)_{t\in [0,T]}$ is a $q$-dimensional Brownian motion. First note
that without loss of generality one may assume that $q=d$ by increasing the dimension of $W$ or adding some identically zero components to $X$ (no ellipticity 
like assumption  is needed here).  This $SDE$ can be
written   in the Stratonovich sense as follows
\begin{equation}\label{strata}
dX_t = f(X_t) \circ dW_t,\quad X_0\!\in L^p_{\R^d}(\P),
\end{equation}
where, for notational convenience $W=(W^0,W^1,\ldots,W^q)$ stands   for $(t,W_t)$, $X_t=(X^0_t,X^1_t, \ldots,
X^d_t)$ stands  for $(t,X_t)$ and $f:[0,T]\times \R^d\to {\cal M}((d+1)\times (d+1))$ (with $f^{0.}(t,x)=(1,0,\ldots,0)$ as $0^{th}$ line) is a differentiable function with bounded differentials. 

Following rough paths theory   initiated by T. Lyons (\cite{LYO}) and developed with many co-authors (see~$e.g.$ \cite{LYO3,  LEJ, LEJ2, LYO3, FRIVIC}  for an introduction), one can also solve this equation
in the  sense of rough paths with finite $p$-variation, $p\ge 2$, since we know ($e.g.$  from the former Kolmogorov criterion) that $W$ $a.s.$ does have 
finite $\frac 1q$-H\" older norm, for any $q>2$.  Namely this means solving an equation formally reading 
\begin{equation}\label{roughpath}
dx_t = f(x_t) d \mathbf{y}_t,\quad x_0\!\in \R^d.
\end{equation}
In   this equation $\mathbf y$ does not represent the path (null at $0$)  $y_t=W_t(\omega)$, $t\!\in [0,T]$
itself but an {\em enhanced  path} embedded  in a larger space, also called {\em geometric multiplicative functional  lying on} $y$  with controlled
$\frac 1q$-H\" older semi-norm, namely a  couple
$\mathbf y=((\mathbf{y}^1_{s,t})_{0\le s\le t\le T},(\mathbf{y}^2_{s,t})_{0\le s\le t\le T}) $ where
$\mathbf{y}^1_{s,t}=y_t-y_s\!\in \R^{d+1},\, 0\le s\le  t\le T$, can be identified with the path $(y_t)$    and
$(\mathbf{y}^2_{s,t})_{0\le s\le t\le T}$ satisfies, $\mathbf{y}^2_{s,t}\!\in \R^{(d+1)^2}$ for every
$0\le s\le u\le t\le T$ and  the following tensor multiplicative property
\[
\mathbf{y}^2_{s,t}=\mathbf{y}^2_{s,u}+\mathbf{y}^2_{u,t}+\mathbf{y}^1_{s,u}\otimes\mathbf{y}^1_{u,t}.
\]  

 Different choices for this
functional are possible, leading to different solutions to the above
Equation~(\ref{roughpath}). The choice that makes  coincide 
$a.s.$ the solution of~(\ref{strata}) and the pathwise solutions of~(\ref{roughpath}) is given by 
\begin{equation}\label{boldy2}
  \mathbf{y}^1_{s,t} = W_t (\omega)-W_s (\omega)
,\; \mathbf{y}^2_{s,t} :=  \left(\int_s^t (W^i_u-W^i_s)\circ dW^j_u\right)_{i,j=0,\ldots, d\hskip -1.25 cm }(\omega)
\end{equation}
so that
\[
\mathbf{y}^1_{s,u}\otimes\mathbf{y}^1_{u,t}= \left(\mathbf{y}^{1,i}_{s,u}\mathbf{y}^{1,j}_{s,u}
\right)_{i,j=0,\ldots, d}.
\]
The term  $\mathbf{y}^2_{s,t}$ is but the ``running" L\'evy areas related to the components of the Brownian mtion $W$. The enhanced path of $W$ will be denoted $\mathbf{W}$ (although we will keep the notation $\mathbf{y}$ in some proofs for notational convenience). One defines, for every $q\ge 1$,  the $\frac 1q$-H\" older distance by setting
\[
\rho_q(\mathbf{y}-\mathbf{x}) = \|\mathbf{y}^1-\mathbf{x}^1\|_{q,Hol}+ \|\mathbf{y}^2-\mathbf{x}^2\|_{q/2,Hol}
\]
where 
\[
\|\mathbf{x}^2\|_{q/2,Hol}:= T^{\frac 2q}\sup_{0\le s< t\le T}\frac{|\mathbf{x}^2_{s,t}|}{|t-s|^{\frac 2q}}.
\]
\ni{\bf Remark.} Likewise, when $p\!\in [2,3)$,  one defines the {\em $p$-variation distance} between two such multiplicative functionals $\mathbf{y}$, $\mathbf{z}$ is defined by
\[
\delta_p(\mathbf{y},\mathbf{z}) ={\rm Var}_{p, [0,T]}(\mathbf{y}^1-\mathbf{z}^1)+{\rm Var}_{p/2,
[0,T]}(\mathbf{y}^2-\mathbf{z}^2) 
\]
where 
$$
\displaystyle  {\rm Var}_{q, [0,T]}(\mathbf{y}^2)\!:=\! 
\sup\!\left\{\!\!\left(\sum_{\ell=0}^{k-1}|\mathbf{y}^2_{t_\ell,t_{\ell+1}}|^q\!\right)^{\frac
1q}\!\!\!\!, 0\!\le\! t_0\!\le \!t_1\!\le\! \cdots \!\le \! t_k\ \le \!T, k\ge 1\!\right\}\!.
$$

The distance $\rho_q$ has been introduced in \cite{LYO95} although rough path theory was originally developed for the distance $\delta_p$ in $p$-variation. Recently several authors came back to H\" older distances $\rho_q$ (see $e.g.$~\cite{LEJ2, FRI, FRIVIC}).

\medskip
The following so-called universal limit theorem  theorem   (including variants)  describes the continuity of the so-called It\^o
map $\mathbf{y}\mapsto x$ with respect to both  $\delta_p$ and $\rho_p$-distances and will be the key for our main result.    It was the starting point of rough path theory initiated by T. Lyons. Several statements (or improvements)    can be found $e.g.$ in~\cite{LYO,  LEJ, LEJb,  LYO3, FRIVIC}. We state here some versions coming from~\cite{LEJ} and \cite{LEJ2}.  

\begin{Theorem} \label{MainThm}   Let $\alpha\!\in (0,1]$.

\smallskip
\noindent $(a)$ (See~\cite{LEJ2}) Let $f: [0,T]\times  \R^d \to {\cal M}((d+1)\times(d+1)$, twice differentiable with a bounded first differential and an $\alpha$-H\" older second differential. Suppose the 
multiplicative functional $\mathbf{y}$  satisfies $ \|\mathbf{y}^1-\mathbf{x}^1\|_{q,Hol}+ \|\mathbf{y}^2-\mathbf{x}^2\|_{q/2,Hol}
 <+\infty$ for  
$q\!\in (2,  2+\alpha)$. Then   Equation~(\ref{roughpath}) has a unique solution starting at $x_0$.

When  $\mathbf{y}= \mathbf{W}(\omega)$ ($i.e.$ given by (\ref{boldy2})),  the first component $\mathbf{x}^1=x$ of the solution solution $\mathbf{x}=(\mathbf{x}^1, \mathbf{x}^2)$ $a.s.$ coincides with $(X_t(\omega))_{t\in [0,T]}$, solution 
to the $SDE$ in the Stratonovich sense. 

Furthermore, the It\^o map $\mathbf{y}\mapsto \mathbf{x}$  is continuous for the  H\" older $\rho_q$ distance (and locally Lipschitz in sense described in~\cite{LEJ2}).

\smallskip
\noindent $(b)$ (See~\cite{FRIVIC, LEJ3})
 If $f\!\in {\cal C}^{2}\big([0,T]\times  \R^d, {\cal M}((d+1)\times(d+1)\big)$  is such that  $f.\nabla f$ is bounded with an $\alpha$-H\"older differential, then the conclusions of claim $(a)$still hold. 
\end{Theorem}

\subsection{Quantization of the $SDE$ and main result}
Let $(\alpha^N)_{N\ge 1}$ denote a sequence of    quantizers of the Brownian motion. Each $\alpha^N$ is made up of $N$ functions (or elementary
quantizer) 
 $\alpha^{N}_n:[0,T]\to
\R^d$, $n=1,\ldots,N$.   For convenience   a component ``$0$" will  be added accordingly  to  each elementary quantizer   
$\alpha_n^{N}$   by  setting $\alpha_n^{N,0}(t)=t$ (which exactly quantizes the function $W^0_t= t$). We assume that every elementary quantizer $\alpha^{N}_n$
is a continuous function with finite variation over $[0,T]$.  The resulting Voronoi quantizer $\widehat W=\widehat W^{\alpha^N}$ of $W$ reads
\[
\widehat W_t = \sum_{n=1}^N \alpha^N_n(t)\mbox{\bf 1}_{\{ W \in C_n(\alpha^{N})\}} ,\quad t\!\in [0,T].
\]

Our aim is to approximate the diffusion process  $(x_t)_{t\in [0,T]}$ solution to the $SDE$~(\ref{strata}) by the solution $\widetilde X^N$ of the equation
\[
d\widetilde X^N_t= f\big(\widetilde X^N_t\big)d\widehat W_t,\; \widetilde X^N_0=x_0.
\]
as $N\to \infty$. In fact, a less formal expression is available for the process $\widetilde X^N$, namely 
\[
\widetilde X^N = \sum_{n=1}^N \widetilde x_n^{N}\mbox{\bf 1}_{\{ W \in C_n(\alpha^{N})\}} 
\]
where each $x^N_n$ is solution to the $ODE$ 
\begin{equation}\label{ODEQuant}
d\widetilde x_n^{N}(t) = f(\widetilde x_n^{N}(t))\,d\,\alpha^{N}_n(t), \quad \widetilde x^{N}_n(0)= x_0,\quad n=1,\ldots,N.
\end{equation}
Note that $X^N$ is  a non-Voronoi quantization of $(x_t)$ (at level $N$). The starting natural idea was to  hope that $X^N$ converges to $(x_t)$ owing to  the convergence of $\widehat W^{N}$ toward $W$\dots in an appropriate sense. Since we will use  the above Theorem~\ref{MainThm}, we need to prove the convergence of the geometric functional $\widehat{\mathbf{W}}^N$ related to $\widehat W$ toward that of $W$. The quantity $\widehat{\mathbf{W}}^N$ is formally defined by mimicking the definition of $\mathbf{W}$, namely, for every $(s,t)\!\in [0,T]$, $0\le s< t\le T$, 
\[
\widehat{\mathbf{W}}^{1,N}\!(\omega) \!:=\!  \widehat W_t (\omega)-\widehat W_s (\omega)
,\; \widehat{\mathbf{W}}^{2,N}_{s,t} \!(\omega)\!:=\!  \left(\int_s^t (\widehat W^i_u-\widehat W^i_s)  d\widehat W^j_u\right)_{i,j=0,\ldots, d}\hskip -1cm (\omega)
\]
still with the convention $\widehat W^{0,N}_t=t$. The integral must be understood in the usual Stieltjes sense.

\begin{Theorem}  \label{MainThm2} Let $(\widehat W^N)_{N\ge 1}$ be a sequence of stationary 
quadratic functional quantizers of the Brownian motion converging to $W$ in $L^2_{L^2_{_T}}(\P)$. 

Let $f$ be like in claims~$(a)$ or $(b)$ in Theorem~\ref{MainThm}. 
Consider for every
$N\ge 1$,  the solutions of the quantized $ODE$
\[
d\widetilde X_t^N = f(\widetilde X^N_t)\, d \widehat W^N_t,\qquad N\ge 1.
\] 
as defined by~(\ref{ODEQuant}). Let $\mathbf{X}$ and $\widetilde{\mathbf{X}}^N$ denote the enhanced paths of $X$, solution to~(\ref{strata}),  and $\widetilde{X}^N$ respectively.    Then, for every $q\!\in(2,2+\alpha)$, 
\[
\rho_q(\widetilde{\mathbf{X}}^N,\mathbf{X} )\stackrel{\P}{\longrightarrow}   0.
\]
Furthermore if $r> \frac{2}{3}$ then
\[
\rho_q(\widetilde{\mathbf{X}}^{\lfloor e^{N^r}\rfloor},\mathbf{X} )\stackrel{a.s.}{\longrightarrow}   0.
\]
\end{Theorem}

In view of what precedes this result is, as announced,  a  straightforward corollary of the continuity of the It\^o map established Theorem~\ref{MainThm}, once the convergence $\rho_q(\widehat{\mathbf{W}}^N,\mathbf{W})$ in probability  is established for any $q\!\in(2,3)$. A slightly more derailed proof is proposed at the end of Section~5.

 In fact we will prove a much precise statement  concerning the Brownian motion  since we will establish for every $q>2$ the convergence in every $L^p(\P)$, $0<p<\infty$, of $\rho_q(\widehat{\mathbf{W}}^N,\mathbf{W})$ with an explicit $L^p(\P)$-rate of convergence in the scale $(\log N)^{-\theta}$, $\theta\!\in (0,1)$. 
 
 These rates can   be transferred to the convergence of the quantized $SDE$, conditionally to some events on which the It\^o map is itself Lipschitz continuous for the distances $\rho_q$.
 Several results of local Lipschitz continuity have been established recently, especially  in~\cite{FRI}, \cite{FRIVIC}, \cite{LEJ2}, \cite{LEJ3}, although not completely satisfactory from a practical point of view. So we decided not to reproduce (and take advantage of) them here.

 The proof is divided into two steps:  the convergence for the H\"older semi-norm) of the regular path component is established in Section~\ref{Kolmopvar} (in which more general processes are considered) and the convergence of approximate L\'evy areas in Section~\ref{Cinq} (entirely devoted to the Brownian case for the sake of simplicity).

\medskip
\noindent {\bf Remarks.} $\bullet$ There is a small abuse of notation in the above Theorem since $\widetilde
X^N$ is not a Voronoi quantizer of $X$: this quantization of $X$ is defined on the Voronoi partition
(for the $L^2_{T,\R^d}$-norm) induced by the quantization of the Brownian motion $W$. 

\medskip
\ni $\bullet$ The same results holds for the Brownian bridge, the Ornstein-Uhlenbeck process and more
generally for continuous Gaussian semi-martingales that satisfy the Kolmogorov criterion.

\section{Convergence  of the paths of processes in H\" older semi-norm}\label{Kolmopvar}

\subsection{A general setting including stationary functional quantization}
In this section we investigate the connections between the
celebrated Kolmogorov criterion and the tightness of some classes
of sequences of processes for the topology of $\frac 1q$-H\"older convergence. In fact this connection is somehow the first step of the rough path theory, but we will look at it in a slightly different way. Whatsoever this naive pathwise convergence   is not sufficient to get the continuity of the It\^o map  in a Brownian framework and we will also have to deal for our purpose with the multiplicative functional (see Section~\ref{Cinq}).

But at this stage we aim at showing that when a sequence $(Y^N)_{N\ge 1}$ satisfies some ``stationarity property" with respect to a process $Y$, several properties of $Y$ can be transferred to the $Y^N$. Indeed, the same phenomenon will occur for the multiplicative function (see the next section).

If  $Y$  satisfies the Kolmogorov criterion and  $({\cal G}_N)_{N\ge 1}$
denotes a sequence of sub-$\sigma$-fields of
${\cal A}$, then a sequence of
processes defined by
\[
 Y^N:=  \E(Y\,|\, {\cal G}^N), \quad N\ge 1,
\]  
where the
conditional expectation is considered in the functional sense (see~Appendix) is ($C$-tight and) tight for a whole
family of  topologies induced by  convergence in $\frac 1q$-H\"older sense.

\begin{Definition} Let  $p\ge 1$, $\theta >0$. A process $Y=(Y_t)_{t\in
[0,T]}$ satisfies the Kolmogorov criterion $ (K_{p,\theta})$ if  there
is a real constant $C^{Kol}_{_T}\!>\!0$ such that 
\[
 \forall\, s,\, t\!\in [0,T], \qquad
\E|Y_t-Y_s|^p\le C^{Kol}_{_T} \, |t-s|^{1+\theta} \qquad \mbox{ and }\qquad Y_0\!\in
L^p(\P).
\]
\end{Definition} 

\begin{Theorem}\label{Thm1}Let $Y:=(Y_t)_{t\in [0,T]}$ be a pathwise continuous
process defined on $(\Omega,{\cal A}, \P)$ satisfying the Kolmogorov
criterion $ (K_{p,\theta})$. Let $({\cal G}_N)_{N\ge 1}$ be a sequence of
sub-$\sigma$-fields of ${\cal A}$. For every $N\ge 1$ set 
\[
Y^N := \E(Y\,|\, {\cal G}_N).
\]

For every $N\ge 1$, $Y^N$ has a pathwise continuous version satisfying
\[
\forall\, t\!\in [0,T],\qquad Y^N_t = \E(Y_t\,|\, {\cal G}_{_N})\quad a.s.
\]
Furthermore, if one of the following conditions is satisfied

\smallskip
\noindent $(a)$ ${\cal G}_N\subset {\cal G}_{N+1}$,

\smallskip
\noindent $(b)$ There exists an everywhere dense  subset $D\subset [0,T]$ such that 
$$
\forall\, t\!\in [0,T], \quad Y^N_t
\stackrel{\P}{\longrightarrow} Y_t.
$$

\smallskip
\noindent $(c)$ $|Y^N-Y|_{L^r_T} \stackrel{\P}{\longrightarrow} 0$ for
some $r\ge 1$,

\medskip
\noindent then
\[
\forall\,q >\frac 1\theta,\qquad \forall\, p\!\in [1, q\theta), \qquad \|Y-Y^N\|_{\rm sup}+\|Y-Y^N\|_{q,Hol} \stackrel{L^p}{\longrightarrow} 0.
\]
\end{Theorem}

The proof of the theorem is a variant of the proof of the
Kolmogorov criterion for functional tightness of processes. It consists in
a string of several lemmas. For the following classical lemma, we refer
to~\cite{LEJ} (where it is stated and proved for semi-norms in $p$-variation).

\begin{Lemma}\label{technicalities} Let $x,\, y\!\in{\cal C}([0,T],\R^d)$
and let $q\ge 1$. Then 

\smallskip
\noindent $(a)$ $\|x-x(0)\|_{\rm sup} \le \|x\|_{q, Hol}$.

\smallskip
\noindent $(b)$ $\|x+y\|_{q, Hol}\le\|x\|_{p, Hol}+\|y\|_{q, Hol}$ if  $q\ge 1$,  

\smallskip
\noindent $(c)$ For every $q>q'\ge 1$, $\|x\|_{q, Hol}^{q} \le (2
\|x\|_{\rm sup})^{q-q'}\|x\|_{q', Hol}^{q'}$. 

\smallskip
\noindent $(d)$ Claims $(a)$-$(b)$-$(c)$  remain true with the $p$-variation semi-norm  ${\rm Var}_{q,[0,T]}$ instead of the $\frac 1q$-H\" older  semi-norm.
\end{Lemma}

\begin{Lemma} \label{Lem2.2} Let $p\!\in [1,\infty)$. If $Y$ satisfies the Kolmogorov
criterion $(K_{p,\theta})$ then, for every $N\ge 1$, the process  $Y^N$ defined by $Y^N_t =\E(Y_t\,|\,{\cal G}_{_N})$ has a continuous modification which is  $\frac{\theta'}{p}$-H\"older   continuous for every $\theta'\!\in (0,\theta)$   ($i.e.$   $\|Y^N\|_{\frac{p}{\theta'}, Hol} <+\infty$ $a.s.$). Furthermore, the sequence $(Y^N)_{N\ge 1}$  is $C$-tight and for every
$\theta'\!\in (0,\theta)$, there exists a random variable
$Z_{\theta'}\!\in L^{p}_{\R}(\P)$ such that 
\begin{equation}\label{Variation1}
  \P(d\omega)\mbox{-}a.s. \; \|Y(\omega)\|_{\frac{p}{\theta'}, Hol} \le Z_{\theta'}
\end{equation} 
and  
\begin{equation}\label{Variation1cond}
\forall\, N\ge 1,\,  \|Y^N(\omega)\|_{\frac{p}{\theta'}, Hol} \le    \E(Z_{\theta'}\,|\, {\cal G}_N)(\omega).
\end{equation} 

In particular, the sequence of H\"older semi-norms $(\|Y^N\|_{\frac{p}{\theta'}, Hol} )_{N\ge1}$  is $L^p$-uniformly integrable.
\end{Lemma}

\noindent {\bf Remark.} As a by-product of the proof we also get that 
\[
\E(Z_{\theta'}^{p} ) \le C_{p,T, \theta,\theta'} C^{Kol}_T
\]
where $C_{T, p,\theta,\theta'}$  is a  finite real constant that only depends upon $p$, $T$, $\theta$ and $\theta'$ (and not on $Y$ or the $\sigma$-fields ${\cal G}_{_N}$).

\bigskip
\noindent {\bf Proof.} First it follows form the Kolmogorov criterion that for every $N\ge 1$, $Y^N$ admits a continuous modification which is $\frac{\theta'}{p}$-H\"older for every $\theta'\!\in (0,\theta)$. Moreover the sequence $(Y^N)_{N\ge 1}$ is $C$-tight since every $Y^N$  satisfies  the same Kolmogorov criterion  $(K_{p,\theta})$ and $Y^N_0= \E (Y_0|{\cal G}_{_N})$ is tight on $\R$ (see~\cite{BIL}, \cite{REYO} p.26). Now, let $s,\,t\!\in [0,T]$, let $m,\, n\ge 1$ be  two  fixed integers. First note that

\begin{eqnarray}
\label{KolmoClassic}\hskip  -2 cm  \sup_{s,t\in [0,T],\, t\le s\le t+ \frac{T}{2^n}}|Y_t-Y_s| &\le & 2\sum_{m\ge 0}\max_{0\le k\le 2^{n+m}-1} |Y_{t^{n+m}_{k+1}}-Y_{t^{n+m}_{k}} |   
 \end{eqnarray}  
 and
\[
\max_{0\le k\le 2^{n+m}-1}|Y_{t^{n+m}_{k+1}}-Y_{t^{n+m}_{k}}|^p \le  \sum_{k= 0}^{2^{n+m}-1} |Y_{t^{n+m}_{k+1}} -Y_{t^{n+m}_{k}}|^p.
\]

For every $\theta'\!\in(0,\theta)$, set
\begin{equation}\label{LeZ}
 Z_{\theta'}:=  \frac 2T  \left(\sum_{n\ge 0}2^{n\frac{\theta'}{p}} \hskip -0.5 cm\sup_{s,t\in [0,T],\, t\le s\le
t+ \frac{T}{2^n}}|Y_t-Y_s|\right). 
\end{equation}
Taking the $L^p$-norm  in~(\ref{KolmoClassic}) yields 
\begin{eqnarray*}
\|Z_{\theta'}\|_{_p} &\le & \left( \frac 2T\right)^{\frac{\theta'}{p}} \sum_{n\ge 0}2^{n\frac{\theta'}{p}} \,\| \sup_{s,t\in [0,T],\, t\le s\le
t+ \frac{T}{2^n}}|Y_t-Y_s|\|_{_p}\\
&\le& 2 \left( \frac 2T\right)^{\frac{\theta'}{p}}\sum_{n\ge 0}2^{n\frac{\theta'}{p}}\sum_{\ell\ge 0}\,\| \max_{0\le k\le 2^{n+m}-1}|Y_{t^{n+m}_{k+1}}-Y_{t^{n+m}_{k}}|\|_{_p}.
\end{eqnarray*} 
On the other hand, owing to the the Kolmogorov criterion $(K_{p,\theta})$, 
\begin{eqnarray*}
\E \max_{0\le k\le 2^{n+m}-1}|Y_{t^{n+m}_{k+1}}-Y_{t^{n+m}_{k}}|^p&\le& \sum_{k=0}^{2^{n+m}-1} \E |Y_{t^{n+m}_{k+1}} -Y_{t^{n+m}_{k}}|^p \\\\
&\le& 2^{n+m} C^{Kol}_T2^{-(n+m)(1+\theta)}T^{-(1+\theta)}\\\\
&=&  C^{Kol}_T  T^{-(1+\theta)}2^{-(n+m)\theta}.
\end{eqnarray*} 
Hence
\[
\E\, Z_{\theta'}^p\le    C^{Kol}_T  C_{p,T,\theta,\theta'} \left(\sum_{n\ge 0} \sum_{m\ge 0}2^{n\frac{\theta'-\theta}{p}}2^{-m \frac{\theta}{p}}\right)^p<+\infty
\]
where the finite real constant $C_{p,T, \theta,\theta'}$ only depends on $p$, $T$, $\theta$ and $\theta'$.
On the other hand, for every $\delta\!\in [0,T]$, there exists a integer $n_{\delta}\ge1$ such that $2^{-(1+n_\delta)}\le
\delta/T
\le 2^{-n_\delta}$. Hence, 
\[
\delta^{-\theta'}\hskip -1 cm \sup_{s,t\in [0,T],\, t\le s\le
t+ \delta}|Y_t-Y_s|^p\le  2^{(1+n_\delta)\theta'} T^{-\theta'}\times \sup_{s,t\in [0,T],\, t\le s\le
t+ \frac{T}{2^n}}|Y_t-Y_s|^p\le 
Z_{\theta'}^{p}.
\]

Consequently, for every $s,\,t\!\in[0,T]$, and every $\omega\!\in\Omega$,
\[
|Y_t(\omega)-Y_s(\omega)|\le Z_{\theta'}(\omega) |t-s|^{\frac{\theta'}{p}}
\] 
$i.e.$
\[
\|Y(\omega)\|_{\frac{p}{\theta'}, Hol}\le Z_{\theta'}(\omega).
\]
Finally, it follows from Jensen's  Inequality that for every $s$, $t\!\in \Q\cap [0,T]$, 
\[
 \P(d\omega)\mbox{-}a.s.\qquad  |Y^N_t(\omega)-Y^N_s(\omega)|\le  \E(Z_{\theta'}\,|\, {\cal G}_N)(\omega)
|t-s|^{\theta'}.
\]
In particular this means that, for every $p\ge 1$ and every $\theta'\!\in (0, \theta)$,
\[
\P(d\omega)\mbox{-}a.s. \qquad 
\|Y^N(\omega)\|_{\frac{p}{\theta'}, Hol}\le   \E(Z_{\theta'}\,|\, {\cal
G}_N)(\omega) <+\infty 
\]
and satisfies the  $L^p$-uniform integrability assumption. $\cqfd$

\bigskip
\noindent {\bf Proof of Theorem~\ref{Thm1}.} The sequence
  $(Y^N)_{N\ge 1}$ being $C$-tight on $({\cal C}([0,T], \R^d), \|\,.\, \|_{\sup})$, so is the case of the 
the sequence
$(Y^N,Y)_{N\ge 1}$  on $({\cal C}([0,T],
\R^{2d}),\|\,.\,\|_{\rm sup})$ since the product topology coincides with the uniform  topology. Let $\Q=
w$-$\lim_{N}\P_{(Y^{N'},Y)}$ denote a weak functional limiting value of
$(Y^N,Y)_{N\ge 1}$. If $\Xi= (\Xi^1,\Xi^2)$ denotes the canonical process on
$({\cal C}([0,T],\R^{2d}),\|\,.\,\|_{\rm sup})$, it is clear that $\Q_{\Xi^2}= \P_{_Y}$. 

\smallskip 
\noindent $\rhd$ {\em Convergence of the sup-norm}.  Assume that $(c)$ holds: the functional
$y\mapsto|y^1(t)-y^2(t)|_{L^r_T}$ is continuous on $({\cal
C}([0,T],\R^{2d}),\|\,.\,\|_{\rm sup})$, consequently, $|\Xi^1-\Xi^2|_{L^r_T}=0$  $\Q$-$a.s.$   $i.e.$ $\Q =
\P_{(Y,Y)}$ so that $(Y^N,Y)\stackrel{{\cal
L(\|\,.\,\|_{\rm sup}) }}{\longrightarrow}  (Y,Y)$ as $N\to \infty$ which simply 
means that
$\|Y^N-Y\|_{\rm sup} \stackrel{\P}{\longrightarrow} 0$. On the other
hand,it follows from Lemma~\ref{Lem2.2} that, for every $N\ge 1$, 
\[
\|Y^N-Y\|^p_{\rm sup} \le C_{p,T} \left(\E(Z_{\theta'}^{p} \,|\, {\cal G}_N)+Z_{\theta'}^{p} \right)\,a.s.
\]
(for a given fixed $\theta'\!\in (0,\theta)$) which implies that $(\|Y^N-Y\|_{\rm sup}^p)_{N\ge 1}$ is uniformly
integrable. Finally,
\[
 \E \,\|Y^N-Y\|^p_{\rm sup} \longrightarrow 0\qquad \mbox{ as }\qquad N\to
\infty.
\]
Assume that $(b)$ holds: it follows that, for every
$t_1,\ldots,t_k\!\in D$, one has \\ $(Y^N_{t_1},\ldots, Y^N_{t_k}) \stackrel{\P}{\longrightarrow}
(Y_{t_1}, \ldots, Y_{t_k})$, which in turn implies that\ the convergence $(Y^N_{t_1},\ldots, Y^N_{t_k}, Y_{t_1}, \ldots, Y_{t_k})  \stackrel{\cal L}{\longrightarrow}
(Y_{t_1}, \ldots, Y_{t_k},Y_{t_1}, \ldots, Y_{t_k})$. This means that $\Q$ and $\P_{(Y,Y)}$ have the same finite dimensional  marginals $i.e.$ $\Q=\P_{(Y,Y)}$. 
One concludes like in $(c)$.

\smallskip 
 If $(a)$ holds, for every $t\!\in[0,T]$, $Y^N_t\to Y_t$
$\P$-$a.s.$, so that $(b)$ is satisfied.

\medskip 
\noindent  $\rhd$ {\em Convergence of the H\" older semi-norm}. Let $q\ge 1$. As concerns the convergence of the $\frac 1q$-H\" older semi-norm, one proceeds as
follows. Let $q'\!\in (\frac{p}{\theta}, q)$ and set $\theta':=
\frac{p}{q'}\!\in(0,\theta)$.  It follows from
Lemma~\ref{technicalities}$(b)$-$(c)$ that
\begin{eqnarray*}
\|Y-Y^N\|_{q,Hol} &\le& 2^{1-\frac{q'}{q}}\|Y-Y^N\|_{\rm
sup}^{1-\frac{q'}{q}}\times \left(\|Y\|_{q',Hol} +\|Y^N\|_{q',Hol} \right)^{\frac{q'}{q}}.
\end{eqnarray*}  
Now let $Z:=Z_{\theta'}$ be defined by~(\ref{LeZ}). 
Then,
\[
\|Y\|_{q',Hol} +\|Y^N\|_{q',Hol}\le  Z+ \big( \E(Z\,|\, {\cal G}_{_N})\big).
\]
 Hence, the sequence $(\|Y\|_{q',Hol} +\|Y^N\|_{q',Hol})_{N\ge 1}$, is tight since it is $L^p$-bounded. On the other hand, $\|Y-Y^N\|_{\rm sup} \stackrel{L^p}{\longrightarrow} 0$ so that \\   $\|Y-Y^N\|_{q,Hol} \stackrel{\P}{\longrightarrow} 0$ as $N\to \infty$. 

Now let $\widetilde \theta =\frac pq\!\in (0,\theta)$. The same argument as above shows that $
\|Y-Y^N\|_{q,Hol}\le \widetilde Z+  \E(\widetilde Z\,|\, {\cal G}_{_N})$  where $\widetilde Z=Z_{\widetilde \theta}$ is still  given by~(\ref{LeZ}). As a consequence, $(\|Y-Y^N\|^p_{q,Hol})_{N\ge 1}$ is uniformly integrable since, for every $N\ge 1$, Jensen's Inequality implies
\[
\|Y-Y^N\|^p_{q,Hol} \le2^{p-1} \big(\widetilde Z^{p}+ \E(\widetilde Z^{p}\,|\, {\cal G}_{_N})\big) 
\]
which finally implies that   $\|Y-Y^N\|_{q,Hol} \stackrel{L^p}{\longrightarrow}0$. $\qquad_\diamondsuit$

\subsection{Application to stationary quantizations of Brownian motion: convergence and rates}

\begin{Theorem}\label{Vit1} $(a)$ Let $(\widehat W^N)_{N\ge 1}$ be a sequence of stationary quadratic functional
quantizers of a standard $d$-dimensional Brownian motion $W$ defined by~(\ref{prodquantifB})
or~(\ref{quantifB}) converging to $W$ in a (purely) quadratic sense, namely $\| \,|W-\widehat W^N|_{L^2_{_T}}\|_2\to 0$ as $N\to \infty$.  Then, for every $q>2$, 
\[
\forall\, p\!\in (0,\infty),\qquad 
\|W-\widehat W^N\|_{q,Hol}\stackrel{L^p}{\longrightarrow}0\; \mbox{ as } \;N\to \infty.
\]

\noindent $(b)$ Let $q>2$. If, for every $N\ge 1$, $\widehat W^N$ is an optimal  {\em product quantization} at level $N$.  Then, for every $p\!\in (0,\infty)$, 
\[
\left\| \|W-\widehat W^N\|_{q,Hol}\right\|_{_p} = o\Big((\log N)^{-\frac 32\min\big(\frac{1}{5}(1-\frac 2q),\frac{1}{p} \big)+\alpha}\Big),\quad \forall\, \alpha>0.
\]
\end{Theorem}

The proof of this Theorem is a consequence of the above Theorem~\ref{Thm1}. So we need to get accurate estimates  for the increments of the processes $W-\widehat W^N$. This is the aim of the following lemma.

\begin{Lemma}\label{LemW1} Let $p\!\in [2,+\infty)$. Let $\widehat W^N$, $N\ge 1$, denote a sequence of optimal {\em product} quadratic quantizers.  For every $\rho\!\in (0,\frac 12)$ and every $\varepsilon \!\in (0,3)$, for every $s,\, t\!\in [0,T],\;s\le t$,
\begin{equation}\label{IneqLpW1}
  \left\| (W_t-W_s)-(\widehat W^N_t-\widehat W^N_s)\right\|_{_p}\le C_{\rho,p,T,d,\varepsilon}|t-s|^{\rho} (\log N)^{-(\frac 12-\rho)\wedge(\frac{3-\varepsilon}{2p})}.
\end{equation}
In particular, if $p\!\in (2,3)$, then
\begin{equation}\label{IneqLpW1b}
  \left\| (W_t-W_s)-(\widehat W^N_t-\widehat W^N_s)\right\|_{_p}\le C_{\rho,p,T,d}|t-s|^{\rho} (\log N)^{-(\frac 12-\rho)}.
\end{equation}
\end{Lemma}

\noindent {\bf Proof.} We may assume without loss of generality that we deal with a one-dimensional Brownian motion $W$, quantized at level $N'=\lfloor N^{\frac 1d}\rfloor$ since everything is done component by component. Set  for every $k\ge 1$, $\widetilde \xi_k := \xi_k -\widehat \xi^{N_k}_k$
where $N_1,\ldots,N_k,\ldots$ denotes the optimal bit allocation of an optimal product quadratic quantization at level $N'$. Keep in mind that for every $k>L_{_{W}}(N')$, $N_k=1$ and that of course $N_1\cdots N_{L_{_{W}}(N')}\le N'$. The random vectors  $(\widetilde \xi_k)_{k\ge 1}$ are independent and centered.

It follows from the $K$-$L$ expansion of $W$ and its product quantization that  
\[
(W_t-W_s)-(\widehat W^{N'}_t-\widehat W^{N'}_s)= \sum_{k\ge 1} \lambda_k \widetilde \xi_k\big(e^W_k(t)-e^W_k(s)\big).
\]
Then, it follows from the B.D.G. Inequality for discrete time martingales that 
\begin{eqnarray*}
  \left\|(W_t-W_s)-(\widehat W^{N'}_t-\widehat W^{N'}_s)\right\|_{_p}&\le &C_{p,T} \left\|\sum_{k\ge 1} \lambda_k \widetilde \xi_k (e^W_k(t)-e^W_k(s))^2\right\|_{\frac p2}^{\frac 12}\\
 &\le& C_{p,T}\left(\sum_{k\ge 1} \lambda_k^{1-\rho}\|\widetilde \xi_k\|_{_p}^2\right)^{\frac 12}|t-s|^{\rho}
\end{eqnarray*}
since, for every $k\ge1$,  
\[
(e^W_k(t)-e^W_k(s))^2 = \frac 8T \sin^2 \Big(\frac{t-s}{\sqrt{\lambda_k}}\Big)\cos^2\Big(\frac{t-s}{\sqrt{\lambda_k}}\Big)\le \frac 8T |t-s|^{2\rho}\lambda_k^{-\rho}.
\]

The random variables $\widehat \xi^{N_k}_k$ being an optimal quadratic quantization of the one-dimensional normal distribution for every $k\!\in \!\{1,\!\ldots\!,L_{_{W}}(N')\!\}$, it  follows from~(\ref{mismatch}) that, there exists for  every $\varepsilon\!\in (0,3)$, a constant $\kappa_{p,\varepsilon}$ such that 
\[
\forall\, m\ge 1,\qquad\| \widetilde \xi_k \|_{_p} =\|\xi-\widehat \xi_k^{N_k}\|_{_p} \le \kappa_{p,\varepsilon} \frac{1}{N_k^{1\wedge \frac{3-\varepsilon}{p}}}
\]
where $\widehat \xi^{m}$ denotes the (unique) optimal quadratic quantization at level $m$ of a normally distributed scalar random variable $\xi$.  As a consequence,
\[ 
\left\|(W_t-W_s)-(\widehat W^{N'}_t-\widehat W^{N'}_s)\right\|_{_p}\le C_{p,T,\varepsilon}|t-s|^{\rho}\left(\sum_{k\ge 1} \lambda_k^{1-\rho}\frac{1}{N_k^{2(1\wedge \frac{3-\varepsilon}{p})}}\right)^{\frac 12}.
\]
 \noindent $\rhd$ Temporarily assume that $p\!\in [2,3)$. One may choose $\varepsilon$ so that 
$1\wedge \frac{3-\varepsilon}{p}=1$.
 Now, keeping in mind that $L':=L_{_{W}}(N')\sim \log N'$ and $\lambda_k\le  c\, k^{-2}$ for a real constant $c>0$,  one gets 
\begin{eqnarray*}
\sum_k \lambda_n^{1-\rho}\frac{1}{N_k^2}&\le & \lambda_{L'}^{\rho}\sum_{k=1}^{L'} \frac{\lambda_k}{N^2_k} + \sum_{k>L'} \lambda_k^{1-\rho}\\
&\le& C_{\rho}\left((\log N')^{2\rho} \sum_{k=1}^{L'} \frac{\lambda_k}{N_k^2} +(\log N)^{2\rho -1}\right).
\end{eqnarray*}
Now,  following $e.g.$~\cite{LUPA1}, we know that the optimal bit allocation yields
\[
\sum_{k=1}^{L'} \frac{\lambda_k}{N_k^2} \le \frac C T(\log N')^{-1}
\]
so that, finally
\[
 \left\|(W_t-W_s)-(\widehat W^{N'}_t-\widehat W^{N'}_s)\right\|_{_p}\le C_{\rho,p,T}|t-s|^{\rho}(\log N')^{\rho -\frac 12}.
 \]
\noindent

\noindent $\rhd$ Assume now that $p\!\in [3,+\infty)$ and $\varepsilon\in (0,3)$.   Set $\tilde p=\frac{p}{3-\varepsilon}>1$ and $\tilde q$ its conjugate exponent. Then, H\" older Inequality implies
 \[
 \sum_{k=1}^{L'} \frac{\lambda^{1-\rho}_k}{N_k^{\frac{2}{\tilde p}} }\le \left(\sum_{k=1}^{L'} \frac{\lambda_k}{N_k^2} \right)^{\frac{1}{\tilde p}}\left(\sum_{k=1}^{L'}\lambda^{1-\frac{\rho p}{p-3+\varepsilon}}_k\right)^{\frac{1}{\tilde q}}.
 \]

We inspect now three possibles cases for $\rho$. 

\smallskip
\noindent $\bullet$  If $0<\rho< \frac 12 (1-\frac{3-\varepsilon}{p})$, then $1-\frac{\rho p}{p-3+\varepsilon}>\frac 12$ so that $\sum_{k\ge 1}\lambda^{1-\frac{\rho p}{p-3+\varepsilon}}_k<+\infty$, which in turn implies that
 \[
  \sum_{k=1}^{L'} \frac{\lambda^{1-\rho}_k}{N_k^{\frac{2}{\tilde p}} }\le C_{\rho,p,T}\Big( \log N'\Big)^{-\frac{3-\varepsilon}{p}}.
 \]
  Furthermore $1-\frac{\rho}{2} >\frac{3-\varepsilon}{p}$.
 
\smallskip
\noindent $\bullet$  If $\frac 12 (1-\frac{3-\varepsilon}{p})<\rho < \frac 12 $, then, $1-\frac{\rho}{2} <\frac{3-\varepsilon}{p}$ and $1-\frac{\rho p}{p-3+\varepsilon}=\frac 12$ so that $\sum_{k\ge1} \lambda^{1-\frac{\rho p}{p-3+\varepsilon}}_k<+\infty$ 
 
 \begin{eqnarray*}
  \sum_{k=1}^{L'} \frac{\lambda^{1-\rho}_k}{N_k^{\frac{2}{\tilde p}} }&\le &C_{\rho,p,T}\Big( \log N'\Big)^{-\frac{3-\varepsilon}{p}}\times \Big (L_{_{W}}(N')^{\frac{2\rho p}{p-3+\varepsilon}-1}\Big)^{1-\frac{3-\varepsilon}{p}}\\
  &=& C_{\rho,p,T}\Big( \log N'\Big)^{2\rho-1}.
\end{eqnarray*}

\noindent $\bullet$  If $\frac 12 (1-\frac{3-\varepsilon}{p})=\rho <\frac 12 $, then $1-\frac{\rho}{2} =\frac{3-\varepsilon}{p}$ and $1-\frac{\rho p}{p-3+\varepsilon}=\frac 12$ so that $\sum_{k=1}^{L'}\lambda^{1-\frac{\rho p}{p-3+\varepsilon}}_k\le C_{\rho,p,T}\log\log N'$ (keep in mind $L'= L_{_W}(N')\sim \log N'$). Hence, for every $\varepsilon'\!\in (0,\varepsilon)$, 
 \[
   \sum_{k=1}^{L'} \frac{\lambda^{1-\rho}_k}{N_k^{\frac{2}{\tilde p}} }= o\Big((\log N')^{-\frac{3-\varepsilon'}{p}}\Big).
 \]
 As  conclusion, we get that
 \begin{eqnarray}
\nonumber \left(\sum_k \lambda_k^{1-\rho}\|\widetilde \xi_k\|_{_p}^2\right)^{\frac 12} &\le &\left (\sum_k \lambda_k^{1-\rho}\frac{1}{N_k^{2(1\wedge \frac{3-\varepsilon}{p})}}\right)^{\frac 12} \\
&= &O\Big( (\log N')^{-(\frac 12-\rho)\wedge(\frac{3-\varepsilon}{2p})}\Big)\label{IneqTech}
 \end{eqnarray}
 which completes the proof since $\log (1+N') \!>\!  \frac 1d \log N$ (which implies $\log N' \!>\!\frac 1d \log(N/2)$).$\cqfd$
 
\bigskip
 \noindent {\bf Proof of Theorem~\ref{Vit1}.} $(a)$  Owing to  the monotonicity of the $L^p$-norms, it is enough to show that, the announced convergence holds for  every $q>2$ and every $p>\frac{2q}{q-2}$ or equivalently for every $p>2$ and every $q>\frac{2p}{p-2}$. This statement follows  for the $\frac 1q$-H\" older (semi-)norm follows from Theorem~\ref{Thm1}$(c)$. Indeed $W$ satisfies the Kolmogorov $K_{p,\theta}$ with $\theta=p/2-1$.  On the other hand, it follows from~\cite{GRLUPA3} that, for any sequence of (Voronoi) quantizations $\widehat W^N$ at level $N$ converging in $L^2_{L^2_{_T}}(\P)$ toward $W$, this convergence also holds in the $a.s.$ sense. So Criterion$(c)$ is fulfilled.

\medskip
\noindent $(b)$  Let $q>2$. The process $W-\widehat W^N $ satisfies $K_{p,\rho p-1}$ for every $\rho \!\in (\frac 1p,\frac 12)$ with ``Kolmogorov constants"
\[
C_{T,p}^{Kol} = C_{p,T,\rho,d, \varepsilon}  (\log N)^{-p[(\frac 12-\rho)\wedge(\frac{3-\varepsilon}{2p})]}, \; \varepsilon \!\in (0,3).
\]
We wish to  apply Lemma~\ref{Lem2.2} (and the remark that follows). 

\smallskip
\noindent $\rhd$  Assume $0<p< \frac{5q}{q-2}$. Then there exists $\eta>0$ such that   $p<p'=\frac{5q}{q-2+\eta}$. Set $\theta'=\frac{p'}{q}$. One cheks that $\frac{1}{p'}+\frac 1q < \frac 12$ so that there exists $\eta'>0$ such that  $\rho=\frac{1}{p'}+\frac 1q +\eta'<\frac 12$. Elementary computations show that $\frac 12-\rho<\frac{3}{2 p}$. Let $\varepsilon\!\in (0,3)$ such that $\frac 12-\rho<\frac{3}{2 p}-\varepsilon$. Consequently,  Lemma~\ref{Lem2.2} (and the remark that follows) imply that 
\[
\left\| \| W-\widehat W^N  \|_{q,Hol}\right\|_{_{p'}}\le C_{q,\eta,\eta',T,\varepsilon}(\log N)^{-(\frac 12-\rho)}
\]
and for  any small enough $\alpha>0$, one my specify $\eta$, $\eta'$ and $\varepsilon$ so that $\frac 12-\rho= \frac{3}{10}(1-\frac 2q)-\alpha$. Finally this bounds holds true for  $p\!\in (0,p')$ since the the $L^p$-norm is non-decreasing.

\smallskip
\noindent $\rhd$ Now, if $p\ge \frac{5q}{q-2}$, one checks that $\frac{3}{2p}\ge\frac 12-\Big( \frac 1p +\frac 1q\Big)$. It becomes impossible to specify $\rho \!\in (0,\frac 12)$ so that $\theta'= \frac pq < \theta =\rho p-1$ and  $1-\rho >\frac {3}{2p}$. So the same specifications as above  lead to 
\[
\left\| \| W-\widehat W^N  \|_{q,Hol}\right\|_{_{p'}}\le C_{q,\eta,\eta', \varepsilon,T}(\log N)^{-\frac {3-\varepsilon}{2p}}
\]
which yields the announced result. $\cqfd$
 
\section{Convergence of stationary quantizations of the Brownian motion for the $\rho_q$-H\"older distance.}\label{Cinq}

In view of what will be needed to apply this theorem to the Brownian motion and its functional quantizations, we need to prove a counterpart of Lemmas~\ref{Lem2.2} and~\ref{LemW1} for $\mathbf{W}^2_{s,t}$. However, for the sake of simplicity, by contrast with the previous section,   we will only deal with the case of the Brownian motion and its stationary quantizations.

The main result of this section is the following Theorem.

\begin{Theorem}\label{Vit2} Let $q>2$.

\smallskip
\noindent  $(a)$ Let $(\widehat W^N)_{N\ge 1}$ be a sequence of stationary quadratic functional
quantizers of a standard $d$-dimensional Brownian motion $W$ defined by~(\ref{prodquantifB})
or~(\ref{quantifB}) converging to $W$ in a (purely) quadratic sense, namely $\| \,|W-\widehat W^N|_{L^2_{_T}}\|_2\to 0$ as $N\to \infty$.  Then,  
\[
\forall\, q>2,\quad \forall\, p>0,\qquad 
 \left\|\rho_q(\mathbf{W},\widehat{\mathbf{W}}^N) \right\|_{_p}\stackrel{}{\longrightarrow}0\quad \mbox{ as } \;N\to \infty.
\]

\noindent $(b)$ Let $q>2$. Assume that, for every $N\ge 1$, $\widehat W^N$ is an optimal  {\em product quantization} at level $N$ of $W$.  
Then, for every $q>2$ and every $p>0$, 
\begin{eqnarray*}
\left\| \|\mathbf{W}^2-\widehat{\mathbf{W}}^{2,N}\|_{\frac q2,Hol}\right\|_{_p}&=& o\Big((\log N)^{-\frac{3}{2}\min\big(\frac{2}{7}(1-\frac 2q),\frac{1}{p} \big)+\alpha}\Big), \;\forall\, \alpha>0,
\end{eqnarray*}
so that, finally,  
\[
\left\|\rho_q(\mathbf{W},\widehat{\mathbf{W}}^N) \right\|_{_p}= o\Big((\log N)^{-\frac 32\min\big(\frac{1}{7}(1-\frac 2q),\frac{1}{p} \big)+\alpha}\Big), \;\forall\, \alpha>0.
\]

\noindent $(c)$  If $r> \frac 23$, then 
\[
\rho_q(\mathbf{W},\widehat{\mathbf{W}}^{\lfloor e^{N^r}\rfloor})= o\Big(N^{-(\frac32r-1)\frac{q-2}{7q}+\alpha}\Big) \;\forall\, \alpha>0,\;\P\mbox{-}a.s.  
\]
\end{Theorem}

Note that the result of interest for our purpose (convergence on multi-dimensional stochastic integrals) corresponds to $q\!\in (2,3)$. The proposition below appears as the counterpart of Lemma~\ref{Lem2.2} on the way to the proof.

\begin{Proposition}\label{Composante2} Let $p>2$. 

\smallskip
\noindent $(a)$ Let $\mathbf{W}^2_{s,t}$ be defined by~(\ref{boldy2}). For every $\tilde \theta'\!\in (0, p-1)$, there exists a random variable $Z_{\tilde \theta'}^{(2)} \!\in L^{p}$ such that 
\[
\P\mbox{-}a.s. \quad \forall\, s,\, t\!\in [0,T], \qquad |\mathbf{W}^2_{s,t}|\le Z_{\tilde \theta'}^{(2)} |t-s|^{\frac{\tilde \theta'}{p}}.
\]

\noindent $(b)$ Let  
$$
\widehat{\mathbf{W}}^{2,N}_{s,t}(\omega)= \left(\int_s^t (\widehat W^i_u-\widehat W^i_s) d\widehat W^j_u\right)_{i,j=0,\ldots, d}\hskip -1.0 cm (\omega),\qquad s,t\,\!\in [0,T], \; s\le t,
$$
where $\widehat W = \widehat W^N$ is a stationary quantization of $W$ (the integration holds in the Stieltjes sense).  Then, for every $p>2$ and every $\tilde \theta'\!\in (0,p-1)$, 
\[
\P\mbox{-}a.s. \quad \forall\, s,\, t\!\in [0,T], \qquad |\widehat{\mathbf{W}}^{2,N}_{s,t}|\le \E(  Z_{p,\tilde \theta'}^{(2)} \,|\,{\cal G}_{_N})|t-s|^{\frac{\tilde \theta'}{p}}. 
\]

\noindent $(c)$ Let $\widetilde{\mathbf{W}}^{2,N}_{s,t}=  \mathbf{W}^2_{s,t}-\widehat{\mathbf{W}}^{2,N}_{s,t}$
where $\widehat W =\widehat W^N$ is now an {\em optimal quadratic product quantization} of $W$ at level $N$. Then, if $p>\frac{1}{\rho}$, for every $\tilde \theta'\!\in (0, p(\rho+\frac 12)-2)$,  for every $\varepsilon \!\in (0,3)$ and every $\delta>0$, there exists a real constant $C_{\rho,p,T,d, \varepsilon, \delta}>0$ such that
\[
\left\|\sup_{s,t\in [0,T]}\frac{|\widetilde{\mathbf{W}}^{2,N}_{s,t}|}{|t-s|^{\frac{\tilde \theta'}{p}}} \right\|_{_p} \le C_{\rho,p,T,d, \varepsilon, \delta} \big(\log N\big)^{-(\frac 12-\rho)\wedge \frac{3-\varepsilon}{2(p+\delta)}}. 
\]
 \end{Proposition}

\noindent {\bf Proof.}  $(a)$ The random variable $Z_{\tilde \theta'}^{(2)}$ of interest is defined by  
\[
Z_{\tilde \theta'}^{(2)} := \frac 2T  \sum_{n\ge 0} 2^{n \frac{\tilde \theta'}{p}} \hskip -0.4cm \sup_{s \le t\le s+\frac{T}{2^n}}|\mathbf{W}^{2}_{s,t}|.
\]
Let $s$, $t\!\in [0,T]$, $s\le t\le s+\frac{T}{2^n}$. We know  from the multiplicative tensor property  that, for every $u\!\in [s,t]$,  
\[
\mathbf{W}^2_{s,t}=\mathbf{W}^2_{s,u}+\mathbf{W}^2_{u,t}+W_{s,u}\otimes W_{u,t} 
\]  
and that,  for every $i$, $j\!\in \{0,\ldots,d\}$,
$$
|W^i_{s,u}\otimes W^j_{u,t} |\le \frac 12(|W^i_{s,u}|^2+| W^j_{u,t} |^2).
$$
To evaluate $\sup_{ t\in [s,s+\frac{T}{2^n}] }\hskip -0.1 cm |\mathbf{W}^{2,}_{s,t}|$, we may restrict to dyadic numbers owing to the continuity in $(s,t)$ of $\mathbf{W}^{2}_{s,t}$. As a consequence, 
we have, still following the classical scheme of Kolmogorov criterion
\begin{eqnarray*}
\sup_{ t\in [s,s+\frac{T}{2^n}] }\hskip -0.25 cm |\mathbf{W}^{2}_{s,t}|&\le &2\sum_{m\ge 0}  \max_{0\le k\le 2^{n+m}-1} |\mathbf{W}^{2}_{t^{n+m}_{k},t^{n+m}_{k+1}}|\\
&&+\max_{0\le k\le 2^{n+m}-1}|W_{t^{n+m}_k,t^{n+m}_{k+1}}|^2.
\end{eqnarray*}
 
 Now
\[
\E\,\max_{0\le k\le 2^{n+m}-1} |\mathbf{W}^{2}_{t^{n+m}_{k},t^{n+m}_{k+1}}|^p\le \sum_{\ell=0}^{2^{m+n}-1} \E\,|\mathbf{W}^{2}_{t^{n+m}_{\ell},t^{n+m}_{\ell+1}}|^p
\]
 and 
 \[
\E\,\max_{0\le k\le 2^{n+m}-1} |W_{t^{n+m}_{k},t^{n+m}_{k+1}}|^p\le \sum_{\ell=0}^{2^{m+n}-1} \E\,|W_{t^{n+m}_{\ell},t^{n+m}_{\ell+1}}|^p
\]
where the norms $|\,.\,|$ are the canonical Euclidean norms on the spaces ${\cal M}((d+1),(d+1))$ and  $\R^{d+1}$ respectively.
 
 It is clear that, for every $i\neq j$, $i$, $j\ge 1$ and every $t\ge s$, 
\begin{eqnarray*}
\|\mathbf{W}^{2,ij}_{s,t}\|_{_p}&=&\left \|\int_s^{t} (W^i_u-W^i_s)dW^j_u\right \|_{_p}\\
&\le &\left \|\int_s^{t} (W^i_u-W^i_s)dW^j_u \right\|_{_p} \\
& \le& C^{BDG}_p \left \|\int_s^{t}(W^i_u-W^i_s)^2du \right\|^{\frac 12}_{_{\frac p2}}\\
&\le &C_p |t'-t|
\end{eqnarray*}
whereas 
$$
\| \,|W_{t'}-W_{t}|^{2}\|_{_p}=  |t'-t|\|\, |W_1|\|_{_p}= C_{p,d} |t'-t|.
$$

 Noting that $W^0_t=t$  and, if $i=j$, $1\le i\le d$,  $\mathbf{W}^{2,ii}_{s,t}= \frac 12 (W^i_t-W^i_s)^2$   shows that the above upper-bound still holds for $i=j$ and $i$ or $j=0$. Consequently, we also have
 \[
 \|\mathbf{W}^{2,ij}_{s,t}\|_{_p}\le C_{p,d} |t'-t|.
 \]
 Consequently
 \[
 \E\!\max_{0\le k\le 2^{n+m}-1} \!|\mathbf{W}^{2}_{t^{n+m}_{k},t^{n+m}_{k+1}}|^p \le C_{p,d}\hskip -0.10 cm  \sum_{k=0}^{2^{n+m}-1} \hskip -0.10 cm \left(\frac{T}{2^{n+m}}\right)^{p}\!=\! C_{p,d,T} 2^{(n+m)(1-p)}
 \]
  so that 
 \[
 \| Z_{\tilde \theta'}^{(2)}\|_{_p} \le C_{p,d,T} \sum_{n\ge 0} 2^{n \frac{\tilde \theta'}{p}}\sum_{m\ge 0}  2^{(n+m)(\frac 1 p-1)} = C_{p,d,T} \sum_{n\ge 0} 2^{n(\frac{\tilde \theta'}{p}-1)}<+\infty 
 \]
since $\tilde \theta'<p-1$.
 
 On the other hand, one has obviously
\[
 \sup_{s,t\in [0,T], s\neq t}\frac{ |\mathbf{W}^{2}_{s,t} |}{|t-s|^{\frac{\tilde \theta'}{p}} } \le  Z_{\tilde \theta'}^{(2)} <+\infty\qquad a.s.
\]
\noindent 
Lemma~\ref{Lem2.2}$(a)$ applied to $W$ (which satisfies $(K_{p, \frac p2-1})$) yields for every $\theta'\!\in (0,\frac p2-1)$ the existence of $Z^{(1)}\!\in L^{p}(\P)$ such that  
\[
 \sup_{s,t\in [0,T], s\le t} \frac{ |\mathbf{W}^{1}_{s,t} |}{|t-s|^{\frac{\theta'}{p}} } \le Z_{\theta'}^{(1)} \qquad a.s.
\]

As a consequence, combining these two results shows that, for every $q>\frac{2p}{p-2}$,
\[
\rho_q(\mathbf{W},0)< Z= Z_{\theta'}^{(1)}+ Z_{\tilde \theta'}^{(2)}\!\in L^p(\P)
\]
where $Z^{(1)}$ is related to $\theta'=\frac pq\!\in (0, \frac p2-1)$ and $Z^{(2)}$ is related to $\tilde \theta'=\frac{2p}{q}\!\in (0,p-2)$.

\smallskip
\noindent $(b)$ If $i\neq j$, $0\le i,j\le d$, it follows from Proposition~\ref{Prop1.1} that $\widehat {\mathbf{W}}^{2,ij,N}_{s,t} = \E(\widehat {\mathbf{W}}^{2,ij,N}_{s,t} \,|\, {\cal G}_N)$ where ${\cal G}_N=\sigma(\widehat W)$ and $\widehat W^N= (\widehat W^{i,N})_{1\le i\le d}$ is an optimal product quantization at level $N$ (which means that for each component $W^i$, $\widehat W^{i,N}$ is an optimal product quantization   at level $N'=\lfloor N^{\frac 1d}\rfloor$).  

\smallskip When $i=j\ge 1$,  $|\widehat{\mathbf{W}}^{2,ii,N}_{s,t}|\le \frac 12 \E\left(  (W^i_t-W^i_s)^2\,|\, {\cal G}_N) \right)$. One derives that
 
\begin{eqnarray*}
\frac{|\widehat{\mathbf{W}}^{2,ii,N}_{s,t}|}{|t-s|^{ \frac{\tilde \theta'}{p} }}& \le & \E \left(\frac{|\mathbf{W}^{2,ii}_{s,t}|}{|t-s|^{\frac{\tilde \theta'}{p} } }\,|\, {\cal G}_N\right)\le    \E\left((Z^{(2)}_{\tilde \theta'})^{\frac{\tilde \theta'}{p}}\,|\, {\cal G}_N\right).
\end{eqnarray*}

When $i=j=0$, $\widehat{\mathbf{W}}^{2,ii,N}= \mathbf{W}^{2,ii}=\frac 12 (t-s)^2$.

\smallskip
\noindent $ (c)$ In this claim, the random variable $Z_{\tilde \theta'}^{(2),N}$ of interest is defined by  
\[
\widetilde Z_{\theta'}^{(2),N}= \frac 2T
\sum_{n\ge 0 }2^{n\frac{\tilde \theta'}{p}} \sup_{s\le t\le s+\frac{T}{2^n}}|\widetilde{\mathbf{W}}^{2,N}_{s,t}|
\]
and we aim at showing that it lies in $L^p(\P)$ with a control on its $L^p$-norm as a function of $N$. One first  derives for $\widetilde{\mathbf{W}}^{2,N}_{s,t}$ the straightforward identity when $s\le u\le t$
\[
\widetilde{\mathbf{W}}^{2,N}_{s,t}=\widetilde{\mathbf{W}}^{2,N}_{s,u}+\widetilde{\mathbf{W}}^{2,N}_{u,t}+\widetilde W_{s,u,t}^N
\]
where 
\begin{eqnarray}
\nonumber \widetilde W^N_{s,u,t}&=&W_{s,u}\otimes W_{u,t}-\widehat W^N_{s,u}\otimes \widehat W^N_{u,t}\\
&=& (W_{s,u}-\widehat W^N_{s,u})\otimes W_{u,t} +\widehat W^N_{s,u}\otimes (W_{u,t}-\widehat W^N_{u,t})\label{identiteWWhat}
\end{eqnarray}
 with $W_{r,s}:=W_r-W_s$ if $r\ge s$, etc. One derives from~(\ref{identiteWWhat}) that
\begin{eqnarray}
|\widetilde{\mathbf{W}}^{2,N}_{s,t}|&\le& 2\sum_{m\ge 0}\max_{0\le k\le 2^{n+m}-1}|\widetilde{\mathbf{W}}^{2,N}_{t^{n+m}_{k},t^{n+m}_{k+1}}|\label{term0}\\
&&\hskip -1cm +2\sum_{m,m'\ge 0}\max_{\substack{0\le k\le 2^{n+m}-1\\0\le k'\le 2^{n+m'}-1}}\hskip-0.25cm |W^{2,N}_{t^{n+m}_{k},t^{n+m}_{k+1}}-\widehat{W}^{2,N}_{t^{n+m}_{k},t^{n+m}_{k+1}}|| W_{t^{n+m}_{k'},t^{n+m}_{k'+1}}|\qquad\label{terme1}\\
&&\hskip -1cm +2\!\!\!\sum_{m,m'\ge 0}\max_{\substack{0\le k\le 2^{n+m}-1\\0\le k'\le 2^{n+m'}-1}   }  \hskip-0.25cm |W^{2,N}_{t^{n+m}_{k},t^{n+m}_{k+1}}-\widehat{W}^{2,N}_{t^{n+m}_{k},t^{n+m}_{k+1}}|| \widehat W_{t^{n+m}_{k'},t^{n+m}_{k'+1}}|.\quad \label{term2}
\end{eqnarray}
where we used that $|u\otimes v|\le |u||v|$. 

We will first deal with deal with the first term in~(\ref{term0}). We note that
\[
\E \,\max_{0\le k\le 2^{n+m}-1}|\widetilde{\mathbf{W}}^{2,N}_{t^{n+m}_{k},t^{n+m}_{k+1}}|^p \le \sum_{0\le k\le 2^{n+m}-1}\E |\widetilde{\mathbf{W}}^{2,N}_{t^{n+m}_{k},t^{n+m}_{k+1}}|^p.
\]

Let $s$, $t\!\in [0,T]$, $s\le t$ and $i$, $j\!\in \{1,\ldots,d\}$, $i\neq j$. One checks that the following decomposition holds
\[
\widetilde{\mathbf{W}}^{2,ij,N}_{s,t}=  \underbrace{\int_{s}^{t} W ^i_{s,u}d(W^j_u-\widehat W ^{j,N}_u)}_{(A)} +  \underbrace{\int_{s}^{t} \widehat W ^{j,N}_{u,t} d(W^{i}_u-\widehat W ^{i,N}_u)}_{(B)}.
\]

Let us focus on $(A)$. First not that, owing to Proposition~\ref{Lem1.1} applied with ${\cal F}_t=\sigma(W^i_u,\, u\in [0,T], \, W^j_s,\, s\le t)$, 
\[
(A) = \sum_{n\ge 1} \widetilde \xi^j_n\int_s^{t} W ^i_{s,u} \cos\Big(\frac{u}{\sqrt{\lambda_n}}\Big)du.
\]
Using that $W^i$ and $W^j$ are independent, one derives that  $(A)$ is the terminal value of   a martingale with respect to the filtration $\sigma(\xi^j_k,\, k\le n,\, W^i_u, 0\le u\le T)$, $n\ge 1$ so that combining    B.D.G. and Minkowski inequalities yields, with the notations of Lemma~\ref{LemW1}, 
\begin{eqnarray*}
\E (|(A)|^p)
&\le& C^{BDG}_{p}  \E \left(\sum_{n\ge 1} (\widetilde \xi_n^j)^2\left( \int_s^{t} W^i_{s,u}\cos\Big(\frac{u}{\sqrt{\lambda_n}}\Big)du\right)^2\right)^{\frac p2} \\
 &\le& C^{BDG}_{p}  \left(\sum_{n\ge 1} \|\widetilde \xi_n^j\|_{_p}^2\left \| \int_s^{t} W^i_{s,u}\cos\Big(\frac{u}{\sqrt{\lambda_n}}\Big)du\right\|_{_p}^2 \right)^{\frac p2}  
\end{eqnarray*}
where $\tilde \xi_n=\xi_n-\widehat \xi_n^{N_n}$ and $N_1,\ldots,N_n,\dots$ denote the optimal bit allocation of an optimal quadratic product quantization at level $N'$ (keep in mind that $N_k=1$, $k>L_{_B}(N')$ and $N_1\cdots N_{L_{_B}(N')}\le N'$ ($B$ scalar Brownian motion). Now an elementary integration by parts yields
\[
 \int_s^{t} W^i_{s,u}\cos\Big(\frac{u}{\sqrt{\lambda_n}}\Big)du= \sqrt{\lambda_n}\int_s^{t}\Big( \sin\Big(\frac{t}{\sqrt{\lambda_n}} \Big)- \sin\Big(\frac{u}{\sqrt{\lambda_n}} \Big) \Big)dW^i_u 
\]
so that, for every $\rho\!\in (0,\frac 12)$, one checks that, owing to the $BDG$ Inequality,  
\[
\left\|   \int_s^{t} W^i_{s,u}\cos\Big(\frac{u}{\sqrt{\lambda_n}}\Big)du \right\|_{_p}\le C^{BDG}_{p}  C_{p,\rho}\lambda_n^{\frac{1-\rho}{2}}|t-s|^{\frac 12+\rho} .
\]

Finally, for every $\varepsilon\!\in (0,3)$,  
\[
\|(A)\|_{_p} \le C_{p,T,\rho,\varepsilon}  \left(\sum_{n\ge 1}  \lambda_n^{1-\rho} \|\tilde \xi_n\|_{_p}^2\right)^{\frac 12} 
|t-s|^{\frac 12+\rho}.
\]
One shows likewise the same inequality  for $(B)$ once noted that
\[
\int_s^{t} \widehat W^{i,N}_{s,u}\cos\Big(\frac{u}{\sqrt{\lambda_n}}\Big)du= \E\left(\int_s^{t}  W^{i}_{s,u}\cos\Big(\frac{u}{\sqrt{\lambda_n}}\Big)du\,|\, {\cal F}^{\widehat W^{i,N}}_T\right)
\]
which implies
\[
\left\|   \int_s^{t}\widehat W^{i,N}_{s,u}\cos\Big(\frac{u}{\sqrt{\lambda_n}}\Big)du \right\|_{_p}\le \left\|   \int_s^{t} W^i_{s,u}\cos\Big(\frac{u}{\sqrt{\lambda_n}}\Big)du \right\|_{_p}.
\]

Consequently, for every $\varepsilon\!\in (0,3)$, 
\begin{eqnarray}
\nonumber  \!\! \!\! \!\!\|\widetilde{\mathbf{W}}^{2,ij,N}_{s,t} \|_{_p}&\le &C_{p,\rho,T}   \left(\sum_{n\ge 1}  \lambda_n^{1-\rho} \|\tilde \xi_n\|_{_p}^2\right)^{\frac 12} 
|t-s|^{\frac 12+\rho}\\
&\le&  C_{p,T,\rho,d,\varepsilon}  \big(\log N\big)^{-(\frac 12-\rho)\wedge \frac{3-\varepsilon}{2p}} 
|t-s|^{\frac 12+\rho}.
\end{eqnarray}

If $i=j\ge1$, then
\[
\widetilde{\mathbf{W}}^{2,ii,N}_{s,t} = \frac12\Big  (\big(W^i_{s,t}\big)^-\big(\widehat W^{i,N}_{s,t})^2\Big)
\]
so that, using again H\" older Inequality, 
\[
\|\widetilde{\mathbf{W}}^{2,ii,N}_{s,t} \|_{_p}= \frac12\|W^i_{s,t}-\widehat W^{i,N}_{s,t}\|_{_{p+\delta)}}\|W^i_{s,t}-\widehat W^{i,N}_{s,t}\|_{_{p(1+\frac{p}{\delta})}}
\]
and one gets the same bounds as in the case $i\neq j$. 

\smallskip
If $i$ or $j=0$, one gets similar bounds: we leave the details to the reader.  Finally, one gets that, for every $i,j\!\in \{0,\ldots,d\}$,

\begin{equation}\label{IneqIncr2}
\nonumber \|\widetilde{\mathbf{W}}^{2,N}_{s,t} \|_{_p}
\le C_{p,\rho,T,d,\varepsilon, \delta} \big(\log N\big)^{-(\frac 12-\rho)\wedge \frac{3-\varepsilon}{2(p+\delta)}}|t-s|^{\frac 12+\rho}.
\end{equation}

\noindent By standard computations similar to those detailed in Lemma~\ref{Lem2.2}, we get
\[
\sum_{m\ge 0} \! \|\!\!\max_{0\le k\le 2^{n+m}\!-1}\hskip -0.30cm|\widetilde{\mathbf{W}}^{2,N}_{t^{n+m}_{k},t^{n+m}_{k+1}}|\|_{_p}\!\le \!C_{p,\rho,T,d,\varepsilon, \delta} \big(\!\log N\!\big)^{-(\frac 12-\rho)\wedge \frac{3-\varepsilon}{2(p+\delta)}}2^{-n(\frac 12 +\rho)}\!.
\]

Let us pass now to the two other sums. We will focus on~(\ref{terme1}) since both behave and can be treated similarly.  
\begin{eqnarray*}
\!\!\max_{\substack{0\le k\le 2^{n+m}-1\\0\le k'\le 2^{n+m'}-1}}\!\!\!|W^{2,N}_{t^{n+m}_{k},t^{n+m}_{k+1}}\!\!&\!-\!&\!\widehat{W}^{2,N}_{t^{n+m}_{k},t^{n+m}_{k+1}}|^p| W_{t^{n+m}_{k'},t^{n+m}_{k'+1}}|^p\\
&\le& \hskip -0.5 cm \sum_{\substack{0\le k\le 2^{n+m}-1\\0\le k'\le 2^{n+m'}-1}}\hskip -0.8 cm |W^{2,N}_{t^{n+m}_{k},t^{n+m}_{k+1}}\!-\!\widehat{W}^{2,N}_{t^{n+m}_{k},t^{n+m}_{k+1}}|^p| W_{t^{n+m}_{k'},t^{n+m}_{k'+1}}|^p\!.
\end{eqnarray*}
Now for every $s,u,t\!\in [0,T]$,  $s\le u\le t$,  it follows from H\" older Inequality that
\begin{eqnarray*}
\|\,|W_{s,u}-\widehat W_{s,u}^N|\,| W_{u,t}|\|_{_p}&\le& \|W_{s,u}-\widehat W_{s,u}^N\|_{_{p+\delta}}\|W_{u,t}\|_{_{p(1+p/\delta)}}\\
&\le& C_{p,\delta} \|W_{s,u}-\widehat W_{s,u}^N\|_{_{p+\delta}}    |t-u|^{\frac 12}.
\end{eqnarray*}

Using   Inequality~(\ref{IneqLpW1}) from Lemma~\ref{LemW1}, we get for every $p>2$, every $\rho\!\in (0,\frac12)$, every $\varepsilon \!\in (0,3)$, 
and every $s,\, t\!\in [0,T]$, $s\le t$, 
\[
\left\| W^i_{s,t}-\widehat W^{i,N}_{s,t}\right\|_{_p}\le C_{\rho,p,T,d,\varepsilon}|t-s|^{\rho} (\log N)^{-(\frac 12-\rho)\wedge(\frac{3-\varepsilon}{2p})}.
\]

Now,
\begin{eqnarray*}
&&\E \hskip -0.4 cm \max_{\substack{0\le k\le 2^{n+m}-1\\0\le k'\le 2^{n+m'}-1}}\hskip -0.25cm |W^{2,N}_{t^{n+m}_{k},t^{n+m}_{k+1}}-\widehat{W}^{2,N}_{t^{n+m}_{k},t^{n+m}_{k+1}}|^p| W_{t^{n+m}_{k'},t^{n+m}_{k'+1}}|^p\\
&&\quad\qquad\le  \Big(C_{\rho,p,T,d,\delta,\varepsilon}  (\log N)^{-(\frac 12-\rho)\wedge(\frac{3-\varepsilon}{2(p+\delta)})}\Big)^p2^{(n+m)(1-\rho p)}2^{(n+m')(1-\frac p2)}
\end{eqnarray*}

\noindent and  we use that  $\rho>\frac 1p$ and $p>2$ to show that
\begin{eqnarray*}
&&\sum_{m,m'\ge 0} \max_{\substack{0\le k\le 2^{n+m}-1\\0\le k'\le 2^{n+m'}-1}}  \| |W^{2,N}_{t^{n+m}_{k},t^{n+m}_{k+1}}-\widehat{W}^{2,N}_{t^{n+m}_{k},t^{n+m}_{k+1}}|^p| W_{t^{n+m}_{k'},t^{n+m}_{k'+1}}|\|_{_p}\\
&&\hskip 4 cm  \le C_{\rho,p,T,d,\delta,\varepsilon}  (\log N)^{-(\frac 12-\rho)\wedge\frac{3-\varepsilon}{2(p+\delta)}}2^{n(\frac2p-(\frac 12+\rho ))}. 
\end{eqnarray*}
 
 Finally, we get
\[
\E \Big(\widetilde Z_{\tilde \theta'}^{(2),N}\Big)^{p}  \le C_{\rho,p,T,d, \varepsilon, \delta} \big(\log N\big)^{-p(\frac 12-\rho)\wedge \frac{3-\varepsilon}{2(p+\delta)}} 
\] 
as soon as $\tilde \theta'\!\in (0,\tilde \theta )$ with $\tilde \theta =p(\rho+\frac 12)-2$. 
Now, it follows by standard arguments  that
\[
\sup_{s,t\in [0,T]}|\widetilde{\mathbf{W}}^{2,N}_{s,t}| \le \widetilde Z_{\tilde \theta'}^{(2),N} |t-s|^{\frac{\tilde \theta'}{p}}
\]
so that, finally
\[
\left\|\sup_{s,t\in [0,T]}\frac{|\widetilde{\mathbf{W}}^{2,N}_{s,t}|}{|t-s|^{\frac{\tilde \theta'}{p}}} \right\|_{_p} \le C_{\rho,p,T,d, \varepsilon, \delta} \big(\log N\big)^{-(\frac 12-\rho)\wedge \frac{3-\varepsilon}{2(p+\delta)}}.   \cqfd
\]

\medskip
Now, we are in position to prove the main result of this section.

\bigskip
\noindent {\bf Proof of Theorem~\ref{Vit2}.} $(a)$  Given Theorem~\ref{Vit1}, this amounts to proving that  $\|\mathbf{W}^2-\widehat{\mathbf{W}}^{2,N}\|_{\frac q2,Hol}$ converges to $0$ in every $L^p(\P)$. This easily follows from Proposition~\ref{Composante2}$(a)$-$(b)$.

\smallskip
\noindent $(b)$  We inspect successively four  cases to maximize $\min (1-\rho, \frac{3}{2p})$ in $\rho$ when it is possible.

\smallskip
\noindent $\rhd$ $q\!\in (2,4)$ and $p<\frac{7q}{2(q-2)}$.  Let $p'$ be defined by $\frac{1}{p'}= \frac{2(q-2)}{7q}+\frac{\alpha}{2}$ with $\alpha>0$ small enough so that  $p'>p$ and $ \frac{1}{p'}+\frac{1}{q}<\frac 12$. Then set $\rho'= \frac{2}{q}+\frac{2}{p'}-\frac 12+\frac{\alpha}{2}$ (note that $\rho'>\frac{1}{p'})$. One checks that $\frac 12-\rho'=1-2(\frac{1}{p'}+\frac{1}{q})= \frac 37(1-\frac 2q)-\alpha\!\in (0,\frac{3-\varepsilon}{2(p'+\delta)}\wedge \frac 12)$ at least for any small enough $\alpha$, $\delta=\delta(\alpha,q)>0$ and $\varepsilon=\varepsilon(\alpha,q)>0$.  Now,  Proposition~\ref{Composante2}$(c)$ applied with $\tilde \theta'= \frac{2p'}{q}< p'(\rho'+\frac 12)-2$   yields the announced asymptotic rate for $\left\| \|\mathbf{W}^2-\widehat{\mathbf{W}}^{2,N}\|_{\frac q2,Hol}\right\|_{_p}$, $p<p'$,  since $L^p(\P)$-norms are non-decreasing in $p$.

\smallskip
\noindent $\rhd$ $q\!\in (2,4)$  and $p\ge \frac{7q}{2(q-2)}$. One sets  the same specifications  as above  for $\rho$ but with $p'=p$. Then  $ 1/2-\rho>\frac{3}{2p}$ and choose $\varepsilon=\varepsilon(q,\alpha)>0$ and $\delta=\delta(q,\alpha)>0$ small enough so that  $ \frac{3-\varepsilon}{2(p+\delta)}\le \frac{3}{2p}+\alpha$.  

\smallskip
\noindent $\rhd$ $q\!\in[4,20/3)$. Then $\frac{7q}{2(q-2)}<\frac{2q}{q-4}$ and one checks that the cases $p\!\in (2,\frac{7q}{2(q-2)})$ and $p\!\in [\frac{7q}{2(q-2)},\frac{2q}{q-4})$ can be solved as above. If $p\ge \frac{2q}{q-4}$ (hence $\ge 5$), no optimization in $\rho$ is possible $i.e.$ any admissible $\rho$ satisfies $\frac 12-\rho>\frac{3}{2p}$.

\smallskip
\noindent $\rhd$ $q\!\ge 20/3$ $i.e.$ $\frac{7q}{2(q-2)}>\frac{2q}{q-4}$.  If $p<\frac{2q}{q-4}$,  set $p'$ such that $\frac{1}{p'} = \frac{q-4}{2q}+\alpha'/2$, $\alpha'>0$ small enough and $\rho'= \frac{2}{q}+\frac{2}{p'}-\frac 12+\frac{\alpha}{2}$. Doing as above yields $\min (1-\rho, \frac{3}{2p})= \frac 2q +\alpha$ for an arbitrary small $\alpha>0$. Note that this quantity is greater than $\frac 37(1-\frac 2q) +\alpha$ (so in that case our exponent is not optimal).  If $p\ge \frac{2q}{q-4}$, we proceed to no optimization in $\rho$.

\smallskip
\noindent $(c)$ This is a  consequence of Borel-Cantelli's Lemma by considering $p>\frac{7q}{q-2}$.$\cqfd$

\bigskip
Now we conclude by proving Theorem~\ref{MainThm2}.

\bigskip
\noindent {\bf Proof of Theorem~\ref{MainThm2}.}  First we check using Proposition~\ref{Composante2} that $\rho_q(\widehat{\mathbf{W}}^N,0)$ and $\rho_q(\mathbf{W},0)$ are $a.s.$ finite since they are integrable. Now we may apply Theorem~\ref{MainThm} which yields the announced result. $\cqfd$
 
\bigskip
\small
\ni {\sc Acknowledgement:} We   thank A. Lejay for helpful discussions and comments about  several versions of this work and   F. Delarue and S. Menozzi for  initiating our first ``meeting"  with rough path theory.

\small
\section*{Appendix: Functional conditional expectation}

Let $(Y_t)_{t\in[0,T]}$ be a bi-measurable process defined on a
probability space  $(\Omega,{\cal A}, \P)$ such that 
\[
\int_0^T \E(Y^2_t) dt <+\infty.
\]
One can consider $Y$ as a   random variable $Y:
(\Omega,{\cal A}, \P)\to L^2_{_T}:= L^2([0,T],dt)$ and more precisely as
an element of the Hilbert space
$$
L^2_{L^2_T}(\Omega,{\cal A},\P):=\left\{Y:(\Omega,{\cal A},
\P)\to L^2_{_T},\; \E \,|Y|_{L^2_T}^2 <+\infty\right\}
$$ 
where $|f|^2_{L^2_T}=\int_0^T f^2(t)dt$. For the sake of simplicity, one 
denotes
$\|Y\|_2:=\sqrt{\E \,|Y|_{L^2_T}^2}$. If ${\cal B}$ denotes a
 sub-$\sigma$-field of ${\cal A}$ (containing all $\P$-negligible sets
of
${\cal A}$) then
$L_{L^2_T}^2(\Omega, {\cal B},\P)$  is a closed sub-space of
$L_{L^2_T}^2(\Omega, {\cal A},\P)$ and one can define the {\em functional
conditional expectation} of 
$Y$ by 
\[
\E(Y\,|\, {\cal B}) := {\rm Proj}_{L_{L^2_T}^2(\Omega,{\cal B},
\P)}^{\perp}(Y).
\]
Functional conditional expectation can be extended to
bi-measurable processes $Y$ such that $\|Y\|_1:=\E\, |Y|_{L^1_T}<+\infty$
following the approach used for $\R^d$-valued random vectors.  Then,
$\E(Y\,|\, {\cal B})$ is characterized by: for every ${\cal B}([0,T])\otimes{\cal B}$-bi-measurable process
$Z=(Z_t)_{t\in[0,T]}$, bounded by $1$, 
\[
 \E \int_0^T Z_t\,Y_t\, dt = \E \int_0^T Z_t \,\E(Y\,|\, {\cal B})_t\,dt.
\] 
In particular, owing to the Fubini theorem,  this implies that as soon as the process
$(\E(Y_t\,|\, {\cal B}))_{t\in [0,T]}$ has a ${\cal B}([0,T])\otimes {\cal B}$ bi-measurable version, the functional conditional expectation could be defined by setting
\[
\E(Y\,|\, {\cal B})_t(\omega) = \E(Y_t\,|\, {\cal B})(\omega),\qquad (\omega,t)\!\in\Omega\times  [0,T].
\]  

\noindent {\sc Examples:} $(a)$ Let ${\cal B} := \sigma({\cal N}_{{\cal
A}},
\, B_i,\; i\!\in I)$ where $(B_i)_{i\in I}$ is a finite measurable
partition of $\Omega$ such that $\P(B_i)>0$, $i\!\in I$.

\smallskip
\noindent $(b)$ Let $Y:=(W_t)_{t\in [0,T]}$ a standard Brownian motion
in $\R^d$ and let ${\cal B}:= \sigma (W_{t_1}, \ldots, W_{t_n})$ where
$0= t_0<t_1<\ldots< t_n =T$. Then
\[
\forall\, t\in [t_k, t_{k+1}), \qquad \E(W\,|\, {\cal B})_t = W_{t_k}
+ \frac{t-t_k}{t_{k+1}-t_k}(W_{t_{k+1}}-W_{t_k}).
\] 


\small\begin{thebibliography}{99.}
\bibitem{BIL} {\sc Billingsley P.} (1968). {\em Convergence of Probability Measure}, Wiley Series in
Probability and Mathematical Statistics, Wiley, New York, 253p. 
\bibitem{COVI}{\sc Coutin L., Victoir N.} (2009).  Enhanced Gaussian Processes and Applications, {\em ESAIM Probab. Stat.}, {\bf  13}: 247--269.
\bibitem{CUMA} {\sc Cuesta-Albertos J.A., Matr{\'a}n C.} (1988). The strong law of large numbers for $k$-means and best
possible nets of Banach valued random variables, {\em Probab. Theory Rel. Fields}, {\bf 78}:523--534.
\bibitem{DER} {\sc Dereich S.} (2008). The coding complexity of diffusion processes under $L^p[0,1]$-norm
distortion, pre-print,  {\em Stoch. proc. and their Appl.}, {\bf 118}(6):938--951.
\bibitem{DEFEMASC} {\sc Dereich S., Fehringer F.,  Matoussi A.,  Scheutzow M.} (2003). On the link between small ball
probabilities and the quantization problem for Gaussian measures on Banach spaces, {\em J. Theoretical
Probab.}, {\bf 16}:249--265.
\bibitem{FRI} {\sc   Friz P.} (2005).  Continuity of the It\^o-map for H\"older rough paths with applications to the support theorem in H\"older norm. {\em Probability and partial differential equations in modern applied mathematics}, 117-135, IMA Vol. Math. Appl., {\bf 140}, Springer, New York.
\bibitem{FRIVIC01} {\sc Friz P.,  Victoir N.} (2007).   Differential Equations Driven by Gaussian Signals I (2007), 
available at arxiv:0707.0313. 
\bibitem{FRIVIC02} {\sc Friz P.,  Victoir N.} (2007).   Differential Equations Driven by Gaussian 
Signals II (2007), available at arxiv:0711.0668. 
\bibitem{FRIVIC} {\sc Friz P.,  Victoir N.} (2010).  {\em Multidimensional Stochastic Differential Equations as Rough Paths: Theory  and Applications}, Cambridge Studies in Advanced Mathematics, 670p.
\bibitem{GRLU} {\sc   Graf S., Luschgy H.} (2000). {\em Foundations of Quantization for Probability
Distributions}. Lect. Notes in Math. 1730, Springer, Berlin.
\bibitem{GRLUPA1} {\sc  Graf S., Luschgy H.,  Pag\`es G.} (2003).  Functional 
quantization and small ball probabilities for Gaussian processes,  {\em J. Theoret. Probab},  {\bf 16}(4):1047--1062.
\bibitem{GRLUPA2} {\sc  Graf S.,  Luschgy H.,  Pag\`es G.} (2007). Optimal quantizers for Radon random vectors in a
Banach space,  {\em J. of Approximation Theory}, {\bf 144}:27--53.
\bibitem{GRLUPA3}  {\sc  Graf S.,  Luschgy H.,   Pag\`es G.} (2008). Distortion mismatch in the
quantization of probability measures, {\em ESAIM P\&S}, {\bf 12}, 127-153.
\bibitem{LEJ} {\sc Lejay A.} (2003). An introduction to rough paths, {\em S\'eminaire de Probabilit\'es YXXVII},
{\em Lecture Notes in Mathematics 1832}:1--59. 
\bibitem{LEJb} {\sc  Lejay A.} (2009). Yet another introduction to rough paths,  {\em S\'eminaire de Probabilit\'es XLII}, {\em Lecture Notes in Mathematics 1979}, 1--101. 
\bibitem{LEJ2} {\sc  Lejay A.} (2009). On rough differential equations,  {\em Electronic Journal of Probability},  {\bf 14}(12):341--364.
\bibitem{LEJ3} {\sc  Lejay A.} (2010). Global solutions to rough differential equations with unbounded vector fields, pre-pub. Inria 00451193, 2010.  
   \bibitem{LUPA1} {\sc Luschgy H.,  Pag\`es G.} (2002). Functional quantization of
stochastic processes, {\em J. Funct. Anal.}, {\bf 196}: 486--531.
\bibitem{LUPA2} {\sc Luschgy H.,  Pag\`es G.} (2004). Sharp asymptotics of the functional quantization
pro\-blem for Gaussian processes, {\em Ann. Probab.} {\bf 32}:1574--1599.
\bibitem{LUPA3}{\sc Luschgy H.,  Pag\`es G.} (2006). Functional
quantization of a class of Brownian diffusions: a constructive approach,  {\em Stoch. Proc. and their
Appl.}, {\bf 116}:310-336.
\bibitem{LUPA4}{\sc Luschgy H.,  Pag\`es G.} (2007).  High resolution product quantization for Gaussian
processes under sup-norm distortion, {\em Bernoulli}, {\bf 13}(3):653--671.
\bibitem{LUPA5}{\sc Luschgy H.,  Pag\`es G.} (2008). Functional
quantization and mean pathwise regularity with an application to L\'evy processes, {\em Annals of Applied Probability}, {\bf 18}(2):427--469.
\bibitem{LUPAWI}{\sc Luschgy H.,  Pag\`es G.,  Wilbertz B.} (2008). Asymptotically optimal quantization schemes for Gaussian processes,  {\em ESAIM P\&S}, {\bf 12}:127--153.
\bibitem{LYO95}{\sc Lyons T.} (1995). Interpretation and Solutions of ODE's Driven by Rough signals, {\em Proc. Symposia Pure}, 1583.
\bibitem{LYO}{\sc Lyons T.} (1998). Differential Equations driven by
rough signals, {\em Rev. Mat. Iberoamericana}, {\bf 14}(2): 215--310.
\bibitem{LYO3} {\sc Lyons T.,  Caruana M.J.,  L\'evy T.} (2007). {\em Differential equations driven by rough paths}.  Lect. Notes in Math. 1908. Notes from T. Lyons's course at \'Ecole d'\'et\'e de Saint-Flour (2004). 
\bibitem{PAPR2}{\sc Pag\`es G.,  Printems J.} (2005).   Functional quantization for numerics with an
application to option pricing,  {\em Monte Carlo Methods \& Applications},  {\bf 11}(4): 407--446.
\bibitem{website} {\sc Pag\`es G.,  Printems J.} (2005). Website devoted to vector and functional optimal
quantization: {\tt www.quantize.maths-fi.com}.
\bibitem{REYO} {\sc  Revuz D.,  Yor M.} (1999). {\em Continuous martingales
and  Brownian motion}. Third edition. Grundlehren der Mathematischen
Wissenschaften [Fundamental Principles of Mathematical Sciences], {\bf
293},  Springer-Verlag, Berlin, 602 p. 
\end{thebibliography}
\end{document}